\documentclass[11pt, hidelinks, oneside]{amsart}

\usepackage[margin=1.3in]{geometry}
\usepackage{amsmath,amsfonts,amssymb,amscd,amsthm,amsbsy}
\usepackage{fancyhdr}
\usepackage{lastpage}
\usepackage{float}
\usepackage{lipsum} 
\usepackage{listings}
\usepackage{caption,subcaption,graphicx,epsfig} 
\usepackage[dvipsnames]{xcolor}
\usepackage{braket}

\definecolor{light-gray}{gray}{0.92}

\usepackage[hidelinks]{hyperref}
\usepackage{tikz-cd,calc}
\usetikzlibrary{arrows.meta}
\usetikzlibrary{knots}
\usetikzlibrary{decorations.markings}

\theoremstyle{plain}
\numberwithin{equation}{section}
\newtheorem{theorem}{Theorem}[section]
\newtheorem{lemma}[theorem]{Lemma}
\newtheorem{proposition}[theorem]{Proposition}

\newtheorem{corollary}[theorem]{Corollary} 

\theoremstyle{definition}
\newtheorem{remark}[theorem]{Remark} 

\newtheorem{example}[theorem]{Example}

\newtheorem{conjecture}[theorem]{Conjecture}

\hypersetup{
colorlinks,
citecolor=blue,
linkcolor=blue
}

\title{Euler class of taut foliations and Dehn filling}

\author[Ying Hu]{Ying Hu}
\address{Department of Mathematics, University of Nebraska Omaha\\
6001 Dodge Street, Omaha, NE 68182-0243}
\email{yinghu@unomaha.edu}
\urladdr{https://yinghu-math.github.io}



\date{\today}
\begin{document}

\maketitle
\begin{abstract}
In this article, we study the Euler class of taut foliations on the Dehn fillings of a $\mathbb{Q}$-homology solid torus. We give a necessary and sufficient condition for the Euler class of a foliation transverse to the core of the filling solid torus to vanish. We apply this condition to taut foliations on Dehn fillings of hyperbolic fibered manifolds and obtain many new left-orderable Dehn filling slopes on these manifolds.  For instance, we show that when $X$ is the exterior of the pretzel knot $P(-2,3,2r+1)$, for $r\geq 3$, $\pi_1(X(\alpha_n))$ is left-orderable for a sequence of positive slopes $\alpha_n$ with $\alpha_0 =2g-2$ and $\alpha_n\to 2g-1$. Lastly, we prove that given any $\mathbb{Q}$-homology solid torus, the set of slopes for which the corresponding Dehn fillings admit a taut foliation transverse to the core with zero Euler class is nowhere dense in $\mathbb{R}\cup \{\frac{1}{0}\}$.
\end{abstract}

\section{Introduction}
\label{section: introduction}
A co-dimension one foliation  on a compact, connected, oriented 3-manifold $M$ is a decomposition of $M$ into a disjoint union of injectively immersed surfaces. These surfaces are called the leaves of the foliation. Throughout the article, we assume that foliations are $C^{\infty, 0}$. So all leaves are smoothly immersed, while the transverse structure may only be $C^0$. By \cite{Cal1}, any topological co-dimension one foliation can be isotoped to be $C^{\infty, 0}$. We also assume that foliations are co-orientable (or transversely orientable). Since $M$ is orientable, this assumption is equivalent to requiring that the tangent plane field of the foliation is orientable.

A closed transversal in a foliated $3$-manifold is a smooth closed loop that is everywhere transverse to the leaves of the foliation. We call a co-dimension one foliation {\it taut} if every leaf of the foliation intersects a closed transversal. See \cite{CKR19} for other notions of tautness for $C^{\infty,0}$ foliations.

\begin{conjecture}[The L-space Conjecture; Conjecture 1 in \cite{BGW}, Conjecture 5 in \cite{Juhasz2015}] 
\label{conj: lspace}
{\it Let $M$ be an irreducible $\mathbb{Q}$-homology sphere. The following statements are equivalent: 

$(1)$ $M$ admits a co-orientable taut foliation.

$(2)$ $\pi_1(M)$ is left-orderable.

$(3)$ $M$ is not an $L$-space.}
\end{conjecture}

A nontrivial group $G$ is called {\it left-orderable} (LO), if there exists a strict total order $<$ on $G$ such that given any elements $a$, $b$ and $c$ in $G$, we have $a<b$ if and only if $ca<cb$. 
 And an {\it L-space} is a $\mathbb{Q}$-homology sphere whose Heegaard Floer homology is of the minimal complexity \cite{OS04threeinvariants,OS04invapplications,OS05lenssurgery,Juhasz2015}. 
 
The conjecture has been confirmed for all graph manifolds \cite{BC, HRRW15, Rasmussen17}. It is also known that $(1) \Rightarrow (3)$ \cite{OS04genusbounds, Bow16, KR17}.  One motivation of our work is to study the equivalence between (1) and (2).  

Given a co-oriented taut foliation $\widehat{\mathcal{F}}$ on an oriented $3$-manifold $M$, the Euler class of the tangent plane field $T\widehat{\mathcal{F}}$, denoted by $e(\widehat{\mathcal{F}})$, is a cohomology class in $H^2(M, \mathbb{Z})$ and the class is zero if and only if the plane field $T\widehat{\mathcal{F}}$ is isomorphic to the product bundle $M\times \mathbb{R}^2\rightarrow M$ (see \S \ref{sec: euler class}). Using Thurston's universal circle action associated with taut foliations (see \cite{CD03}), we have the following theorem whose proof can be found in \cite[\S 7]{BH19}.

\begin{theorem}[Thurston; Calegari-Dunfield, Boyer-Hu]
Let $M$ be a $\mathbb{Q}$-homology $3$-sphere. Suppose that $M$ admits a co-oriented taut foliation whose Euler class is zero. Then $\pi_1(M)$ is left-orderable. 
\label{thm: universal circle}
\end{theorem}

 Theorem \ref{thm: universal circle} leads us to the natural question: when is the Euler class of a given taut foliation on a $\mathbb{Q}$-homology sphere zero?

\subsection{When is the Euler class zero?}
We first note that the Euler class of any oriented tangent plane field over a connected compact $3$-manifold is an ``even'' class (see Proposition \ref{prop:Euler class is even}). As a result,  if $H^2(M)$ is a direct sum of $\mathbb{Z}_2$, then the Euler class of any taut foliation on $M$ must be zero.

\begin{corollary}
Let $M$ be an oriented $3$-manifold satisfying that $H^2(M)$ is isomorphic to a (possibly trivial) direct sum  of $\mathbb{Z}_2$. If there exists a co-orientable taut foliation on $M$, then $\pi_1(M)$ is left-orderable.
\label{cor: euler class zero z2 torsion}
\end{corollary}

In this article, we view $\mathbb{Q}$-homology spheres as Dehn fillings on $\mathbb{Q}$-homology solid tori. See \S \ref{subsec: notation and orientation} for the convention of our notations. For simplicity, we will only state the result for Dehn fillings of $\mathbb{Z}$-homology solid tori below.  We refer the readers to Theorem \ref{thm: e=0 necessary sufficient condition} and Remark \ref{rem: e(F) = 0 iff condition k>1} for more general statements.  

\begin{theorem}
Let $X$ be a $\mathbb{Z}$-homology solid torus and $F$ be a properly embedded surface in $X$ representing a generator of  $H_2(X,\partial X)$. Suppose that $\widehat{\mathcal{F}}$ is an oriented co-dimension one foliation on $X(p/q)$, $p>0$ whose restriction  to the filling solid torus $N$ is  the foliation by meridian disks, and the orientation of the leaves of $\widehat{\mathcal{F}}$ agree with the given orientation of the meridian disk of $N$. Let $\mathcal{F} = \widehat{\mathcal{F}}|_X$ and $\sigma$ denote a nowhere vanishing outward pointing section of $T\mathcal{F}$ along $\partial X$.  Then the Euler class $e(T\widehat{\mathcal{F}}) = 0$ in $H^2(X(p/q))$ if and only if $aq  \equiv 1 \mbox{ \rm{ (mod $p$)}}$, where $e_\sigma(\mathcal{F}) \in H^2(X,\partial X)$ is the relative Euler class of $\mathcal{F}$ associated with the section $\sigma$ and $a = e_\sigma(\mathcal{F})([F])$.
\label{thm: euler class zero knots in zhs}
\end{theorem}

The condition that the restriction of $\widehat{\mathcal{F}}$ to the filling solid torus $N$ is the foliation by meridian disks is equivalent to that the core of $N$ is transverse to $\widehat{\mathcal{F}}$ by shrinking $N$ if necessary. Throughout the article, we always assume that $\widehat{\mathcal{F}}|_{N}$ is the product foliation  if $\widehat{\mathcal{F}}$ on $X(\alpha)$ is transverse to the core of $N$. The existence of taut foliations satisfying this transverseness condition in $X(p/q)$ is also equivalent to that the slope $p/q$ is {\it strongly CTF-detected} as defined in \cite[Definition 6.5]{BC}. 

Almost all known taut foliations on Dehn filled manifolds are transverse to the core of the filling solid torus, and hence Theorem \ref{thm: euler class zero knots in zhs} applies. These include taut foliations constructed in \cite{Roberts1995,RobertsSurfacebundle1, RobertsSurfacebundle2, LR, krishna18, DRcompositeknots19, DRdoublediamond19}, among which the relative Euler classes $e_\sigma(\mathcal{F})$  in \cite{Roberts1995,RobertsSurfacebundle1, RobertsSurfacebundle2, LR, krishna18} are of Thurston norm one.

\subsection{Dehn fillings of fibered manifolds and left-orders}
In \cite{RobertsSurfacebundle1, RobertsSurfacebundle2, krishna18}, families of taut foliations are constructed on manifolds obtained from Dehn fillings of the exteriors of fibered knots.  In the theorem below, the {\it canonical meridian $\mu_0$} and the {\it degeneracy slope} are defined in \S \ref{subsec: fibered knot}. When $X$ is the exterior of a fibered knot in $S^3$, $\mu_0$ is mostly the meridian of the knot \cite[Proposition 7.4]{RobertsSurfacebundle2}. In particular, this is the case when $K$ is the closure of a positive braid considered in Theorem \ref{thm: foliation 3 braids}. 
 
\begin{theorem}[Theorem 4.7 in \cite{RobertsSurfacebundle2}]
 Let $X$ be the exterior of a hyperbolic fibered knot in a closed $3$-manifold $M$ with degeneracy slope $\gamma$ and $\mu_0$ be the canonical meridian of $X$.
 \begin{enumerate}
  \item If $\gamma = \mu_0$, then there exists a taut foliation on $X(\alpha)$ transverse to the core of the filling solid torus for any slope $\alpha \neq \mu_0$.  
  \item If $\gamma$ is a positive slope with respect to $\mu_0$, then there exists a taut foliation on $X(\alpha)$ transverse to the core of the filling solid torus for any $\alpha \in (-\infty, 1)$. 
  \item If $\gamma$ is a negative slope with respect to $\mu_0$, then there exists a taut foliation on $X(\alpha)$ transverse to the core of the filling solid torus for any $\alpha \in (-1, \infty)$.  
 \end{enumerate}
 \label{thm: foliation fibered knot}
\end{theorem}

\begin{theorem}[Theorem 1.2 in \cite{krishna18}]
Let $X$ be the exterior of the closure of a positive  $3$-braid in $S^3$. There exists a taut foliation on $X(\alpha)$ transverse to the core of the filling solid torus  for every slope $\alpha < 2g - 1$, where $g$ is the genus of the closed braid. 
\label{thm: foliation 3 braids}
\end{theorem}

Taut foliations in the theorems above are carried by very nice branched surfaces (\S \ref{subsec: euler class branched surfaces}), which allows us to compute their Euler classes precisely.

\begin{theorem}
Suppose that $\widehat{\mathcal{F}}$ is a taut foliation given by Theorem \ref{thm: foliation fibered knot} or Theorem \ref{thm: foliation 3 braids}. Let $\alpha = p\mu_0 + q\lambda$, where $p>0$ and $\mu_0$ be the canonical meridian. Then $e(\widehat{\mathcal{F}}) = 0$ in $H^2(X(\alpha))$ if and only if  $e(\mathcal{F}) =0$ in $H^2(X)$ and
\begin{equation}
\label{equ: fiber knot slope equation}
 (2g-1)|q| \equiv 1 \mbox{ \rm{(mod $p$)}},
\end{equation}
where $g$ is the genus of the knot and $\mathcal{F} = \widehat{\mathcal{F}}|_X$ is the restriction of $\widehat{\mathcal{F}}$ to the knot exterior $X$.
\label{thm: euler class Rachel's foliation}
 \end{theorem}

For each $g>0$, we denote the set of slopes that satisfies Equation (\ref{equ: fiber knot slope equation}) in Theorem \ref{thm: euler class Rachel's foliation} by 
\begin{displaymath}
\mathcal{L}_g = \{p/q \in \mathbb{Q}\setminus \{0\} : (2g-1) |q|  \equiv 1 \mbox{ \rm{(mod $p$)}}\}.
\end{displaymath}
 We collect some key properties of the set $\mathcal{L}_g$ below. See \S \ref{sec: computation of the slopes} for more details.

\begin{enumerate}
\item The set $\mathcal{L}_1$ contains all nonzero integer slopes. In fact, a slope is in $\mathcal{L}_1$ if and only if it is a nonzero integer slope with respect to some meridian of $X$.  (Lemma \ref{lem: slope branched surface g=1}).
\item When $g>1$, the set $\mathcal{L}_g$ is bounded between $\pm (2g-1)$. There are two sequences of slopes $\pm \{\alpha_n\}$ in  $\mathcal{L}_g$ satisfying $\alpha_0 = 2g-2$ and $\alpha_n \to (2g-1)$ as $n\to \infty$ (Lemma \ref{lem: slope branched surface g>1}).
\end{enumerate}

From Theorem \ref{thm: universal circle}, we obtain the following results on the left-orderability of $3$-manifold groups. 

\begin{corollary}
 Let $K$ be a hyperbolic genus one fibered knot in a $\mathbb{Q}$-homology sphere with the degeneracy slope $\gamma$. With respect to the canonical meridian $\mu_0$,  we have
 \begin{enumerate}
  \item If $\gamma = \mu_0$, then given any $\alpha \in \mathcal{L}_1$, there exists a taut foliation on $X(\alpha)$ of zero Euler class and $\pi_1(X(\alpha))$ is LO.  In particular, $\pi_1(X(\alpha))$ is LO for all integer slopes. 
  \item If $\gamma$ is a positive slope, then given any $\alpha \in \mathcal{L}_1 \cap (-\infty, 1)$, there exists a  taut foliation on $X(\alpha)$ of zero Euler class and $\pi_1(X(\alpha))$ is LO. In particular, $\pi_1(X(\alpha))$ is LO for all non-positive integer slopes. 
  \item If $\gamma$ is a negative slope, then given any $\alpha \in \mathcal{L}_1 \cap (-1, \infty)$, there exists a  taut foliation on $X(\alpha)$ of zero Euler class and $\pi_1(X(\alpha))$ is LO.  In particular, $\pi_1(X(\alpha))$ is LO for all non-negative integer slopes. 
 \end{enumerate}
 \label{cor: left-orderable group g=1 fibered}
\end{corollary}

\begin{remark}[Group left-orderability and genus one fibered knots]
\label{rem: lo and genus one}
Roberts and Shareshian have shown that $\pi_1(X(\alpha))$ is not LO when $\alpha\in [1,\infty)$ in Case (2) \cite[Theorem 1.3]{RS10}, and therefore the analogous result holds in Case (3) (see Remark \ref{rem: genus one trace}).  Though the L-space Conjecture predicts that $\pi_1(X(\alpha))$ is LO whenever $\alpha\in (-\infty, 1)$  in Case (2) and $\alpha\in (-1, \infty)$ in Case (3), this was previously unknown even for integer slopes. 

In Case (1),  Fenley first showed that $\pi_1(X(\alpha))$ is LO for any $\alpha\in \mathbb{Z}$ by showing the existence of $\mathbb{R}$-covered Anosov flows on these manifolds (\cite[Theorem D]{Fen94}; also see \cite[Proposition 3.3]{IN2020}).  Recently, Zung has significantly generalized Fenley's approach to pseudo-Anosov flows and non-$\mathbb{R}$-covered taut foliations. It follows from \cite[Theorem 1]{Zung2020} that $\pi_1(X(\alpha))$ is LO for any slope $\alpha\neq \mu_0$ when $\gamma = \mu_0$.  
\end{remark}

\begin{corollary}
Let $X$ be the exterior of the closure of a positive  $3$-braid in $S^3$ of genus $g>1$. For any slope $\alpha\in \mathcal{L}_g\cap (-\infty, 2g-1)$, there exists a taut foliation on $X(\alpha)$ of zero Euler class and hence $\pi_1(X(\alpha))$ is LO. In particular, there is a monotone increasing sequence of positive slopes $\{\alpha_n\}$ with $\alpha_0 = 2g-2$ and $\alpha_n\to 2g-1$ as $n\to \infty$, satisfying $\pi_1(X(\alpha_n))$ is LO for any $n\in \mathbb{N}$.
\label{cor: 3 braids left-orderability}
\end{corollary}

\begin{remark}[Group left-orderability and L-space knots] 
A positive L-space knot is a knot in $S^3$ whose exterior admits a positive L-space Dehn filling \cite{OS05lenssurgery}.  It is known that the set of finite L-space filling slopes for a genus $g$ positive L-space knot is precisely $[2g-1, \infty)$ (\cite{OS05lenssurgery, KMOS07},\cite[Theorem 1.6]{RR17}). According to the L-space Conjecture, one expects  that given $X$ the exterior of a positive L-space knot, the closed manifold $X(\alpha)$ admits taut foliations and $\pi_1(X(\alpha))$ is LO for any slope $\alpha \in (-\infty, 2g-1)$. 

The closures of positive $3$-braids are examples of positive L-space knots, among which the pretzel knots $P(-2,3,2r+1)$, $r\geq 3$ are the most studied and have been recognized as the only hyperbolic positive L-space Montesinos knots \cite{LM16,BM18}.  In Theorem \ref{thm: foliation 3 braids}, Krishna constructed a taut foliation on $X(\alpha)$ for any $\alpha\in (-\infty, 2g-1)$. Moreover, Nie has shown that $\pi_1(X(\alpha))$ is not LO for any $\alpha\in [2g-1, \infty)$ \cite[Theorem 2]{Nie19}. So to confirm the L-space conjecture for manifolds obtained from Dehn surgeries along pretzel knots $P(-2,3,2r+1)$, $r\geq 3$, it remains to show that $\pi_1(X(\alpha))$ is LO for any $\alpha\in (-\infty, 2g-1)$. However, in general, this was confirmed only for slopes $\alpha$ in an open neighborhood of $0$ (\cite[Theorem 3]{Nie19}; \cite{CD18, HZ19}), except for the case when $r=3$. 

When $X$ is the exterior of the pretzel knot $P(-2,3,7)$, Culler and Dunfield proved that $\pi_1(X(\alpha))$ is LO for any $\alpha\in (-\infty, 6)$ \cite[Figure 3]{CD18}. The genus of $P(-2,3,2r+1)$ is $r+2$, so for $r=3$, the expected LO filling slope interval is  $(-\infty, 9)$. Corollary \ref{cor: 3 braids left-orderability} shows that  there is an increasing sequence of slopes $\alpha_n$ with $\alpha_0 = 8$ and $\alpha_n\to 9$ such that $\pi_1(X(\alpha_n))$ are LO. 
\end{remark}

\subsection{Restrictions on the filling slopes}
Given a $\mathbb{Q}$-homology solid torus $X$ ,  we define $\mathcal{S}_X$ to be the set of slopes satisfying: given any $
\alpha\in\mathcal{S}_X$, there exists a taut foliation on $X(\alpha)$ transverse to the core of the Dehn filling solid torus whose Euler class is zero.

By fixing a meridian $\mu$, we identify the set $\mathcal{S}_X$ with a subset of $\mathbb{Q}\cup\{\frac{1}{0}\}$, denoted by  $\mathcal{S}_{X, \mu} \subset\mathbb{Q}\cup\{\frac{1}{0}\}$. For simplicity, we state the result when the longitude of $X$ is null-homologous in $H_1(X)$. See Theorem \ref{thm: restriction on SX} for the case when $X$ is an arbitrary $\mathbb{Q}$-homology solid torus. 

\begin{theorem}{\label{thm: restriction on SX null homologous}}
Let $X$ be the exterior of a null-homologous knot of genus $g>0$ and $\mu$ be any fixed meridian. 
\begin{enumerate}
 \item Outside the interval $(-2g, 2g)$, the set $\mathcal{S}_{X, \mu}$ contains only $\mu$ and the integer slopes. That is, 
 $$\mathcal{S}_{X, \mu} \setminus (-2g, 2g)\subseteq \mathbb{Z}\cup \{\mu \}.$$
 \item The set $\mathcal{S}_{X,\mu}$ is nowhere dense in $\mathbb{R}\cup\{\frac{1}{0}\}\cong \mathbb{R}P^1$. Particularly, it is nowhere dense in $(-2g, 2g)$.
\end{enumerate}
\end{theorem}

When we use a different meridian, a slope $\frac{p}{q}$ becomes $\frac{p}{q+np}$ for some $n\in \mathbb{Z}$. Since Theorem \ref{thm: restriction on SX null homologous} holds regardless of the choice of the meridian, one can obtain a stronger version of Theorem \ref{thm: restriction on SX null homologous}(1) (see Corollary \ref{cor: obstruction k>1}).

\subsection*{Organization of the paper.}
In \S \ref{sec: euler class}, we review the definition of  the (relative) Euler class and prove Corollary \ref{cor: euler class zero z2 torsion}. Section \ref{sec: euler class and dehn filling} contains the core of our computations (Lemma \ref{lem:cohomology class zero}).  Theorem \ref{thm: e=0 necessary sufficient condition} is also proven in this section. In \S \ref{sec: examples}, we apply the results in \S \ref{sec: euler class and dehn filling} to studying the Euler classes of taut foliations on Dehn fillings of fibered manifolds. We prove Theorem \ref{thm: euler class Rachel's foliation} as well as its corollaries on the left-orderability of $3$-manifold groups in this section.   Finally, Theorem \ref{thm: restriction on SX null homologous} and its generalization (Theorem \ref{thm: restriction on SX}) are proved in \S \ref{sec: euler class slopes}.

\section{The Euler class of tangent plane fields over $3$-manifolds}
\label{sec: euler class}
\subsection{The (Relative) Euler class}
\label{subsec: define euler class}
Let $\xi : E\to B$ be an oriented $2$-plane vector bundle over an oriented CW complex $B$. A $2$-plane bundle is called trivial if it is isomorphic to the product bundle $B\times \mathbb{R}^2\rightarrow B$. It is easy to see that an oriented $2$-plane bundle $\xi$ is trivial if and only if there exists a nowhere vanishing section $\sigma: B\rightarrow E$. In fact, given such a section $\sigma$, it determines an orientable line bundle $\gamma$. Then $\xi \cong\gamma\oplus \xi/\gamma$ is isomorphic to the sum of two orientable line bundles, and hence it is trivial.

We denote the Euler class of an oriented $2$-plane bundle $\xi: E\rightarrow B$ by $e(\xi)$.  A representative cocycle of the cohomology class $e(\xi)$ can be constructed as follows. 

First note that since $\mathbb{R}^2\setminus \{0\}$ is a $K(\mathbb{Z},1)$ space, the only obstruction to the existence of a nowhere vanishing section $\sigma$ arises when one tries to extend a section over the $1$-skeleton $B^{(1)}$ of $B$ to the $2$-skeleton $B^{(2)}$. Fixing a nowhere vanishing section $\sigma: B^{(1)}\rightarrow E$ over the $1$-skeleton $B^{(1)}$, we construct  a cellular $2$-cochain $c_\sigma: C_2(B)\rightarrow \mathbb{Z}$ by specifying its value on  each $2$-cell. Let  $\varphi_\alpha: D_\alpha^2\rightarrow B$ be the characteristic map of a $2$-cell $e_\alpha$ and let $E_{D_\alpha^2}\rightarrow D_\alpha^2$ denote the pullback of $E\rightarrow B$ through $\varphi_\alpha$.  Since $D_\alpha^2$ is contractible, $E_{D_\alpha^2}\rightarrow D_\alpha^2$ is trivial. We identify $E_{D_\alpha^2}$ with the product $D_\alpha^2\times \mathbb{R}^2$ by fixing a trivialization. The nowhere vanishing section $\sigma: B^{(1)}\rightarrow E$ determines a nowhere vanishing section of $E_{D_\alpha^2}\rightarrow D_\alpha^2$ along $\partial D_\alpha^2$, which we also denote by $\sigma$. Since $\sigma$ is nowhere vanishing, the image of the composition $\partial D_\alpha^2\xrightarrow{\,\sigma\,} E_{D_\alpha^2}\xrightarrow{\,\cong\,} D_\alpha^2\times \mathbb{R}^2 \rightarrow \mathbb{R}^2$ is contained in $\mathbb{R}^2\setminus 0$. We set the value of the $2$-cochain $c_\sigma$ on the $2$-cell $e_\alpha$ to be the degree of the composite map $\partial D_\alpha^2\rightarrow \mathbb{R}^2\setminus \{0\}$.  See the diagram below. 

\begin{center}
\begin{tikzpicture}[scale=0.88]
% \draw [help lines] (0,0) grid (12, 6);
 \node at (2,4) {$\mathbb{R}^2\setminus 0$};
 \draw [thick, ->] (3.8, 4) -- (2.8, 4);
 \node at (4.7, 4) {$D^2_\alpha\times \mathbb{R}^2$};
 \draw [thick, ->] (6.6, 4) -- (5.6,4);
 \node [above] at (6.1,4) {$\cong$};
 \node at (7.1,4) {$E_{D^2_\alpha}$};
 \draw [thick, ->] (7.5, 4)--(8.5,4);
 \node at (8.8,4) {$E$};
 \draw [thick, ->] (8.8, 3.7) -- (8.8, 2.5);
 \node at (8.8, 2.2) {$B$};
 \draw [thick, ->] (7.5, 2.2)--(8.5,2.2);
 \node [above] at (8, 2.2) {$\varphi_\alpha$};
 \node at (7.1, 2.2) {$D_\alpha^2$};
 \draw [thick, ->] (7.1, 3.7) -- (7.1, 2.5);
 \node at (5.1, 2.2) {$\partial D_\alpha^2$};
 \draw [thick, ->] (5.6, 2.5)--(6.8, 3.7);
 \node [above] at (6, 3) {$\sigma$};
  \draw [thick, <-] (6.7, 2.2) -- (5.7,2.2);
\end{tikzpicture}
\end{center}

The $2$-cochain $c_\sigma$ is a cocycle and the Euler class $e(\xi)$ is defined to be  the cohomology class $[c_\sigma]$ in $H^2(B)$. See \cite[Section 12]{Mil74} or \cite[Chapter 4]{CC03} for instance. 

Given  a subcomplex $A$ of $B$, suppose that $\xi: E\rightarrow B$ is trivial over $A$. So there exists a nontrivial section of $\xi$ over $A$. Then if one starts with a section $\sigma: B^{(1)}\cup A\rightarrow E$, the $2$-cocyle $c_\sigma$ constructed above vanishes over $A$ and hence defines a relative cohomology class $[c_\sigma]$ in $H^2(B,A)$. We call the relative class $[c_\sigma]$ the relative Euler class of $\xi$ associated to the section $\sigma$, denoted by $e_\sigma(\xi)$. The relative class vanishes if and only if there exists a section of $\xi$ that extends $\sigma$. In general, $e_\sigma(\xi)$ does depend on  $\sigma$. By construction, the image of $e_\sigma(\xi)$ under the homomorphism $H^2(B,A)\rightarrow H^2(B)$ is the Euler class $e(\xi)$ in $H^2(B)$.

\subsection{Euler class and $\mathbb{Z}_2$-torsion}
\label{subsec: euler class Z2 torsion}
Proposition \ref{prop:Euler class is even} below  states a general fact of the Euler class of an oriented tangent plane field over an oriented $3$-manifold. Note that Corollary \ref{cor: euler class zero z2 torsion} follows immediately from Proposition \ref{prop:Euler class is even} and Theorem \ref{thm: universal circle}.

\begin{proposition}
  Let $M$ be an oriented closed $3$-manifold and $\xi: E\rightarrow M$ be an oriented tangent plane field of $M$. Then the Euler class $e(\xi)$ lies in the image of the homomorphism $H^2(M,\mathbb{Z})\xrightarrow{\cdot 2} H^2(M,\mathbb{Z})$ induced by $\mathbb{Z}\xrightarrow{\cdot 2} \mathbb{Z}$. In particular, if $H^2(M)$ is isomorphic to a direct sum  of $\mathbb{Z}_2$, then $e(\xi)=0$. 
  \label{prop:Euler class is even}
\end{proposition}

Before we prove the proposition, we briefly review  the Stiefel-Whitney classes. We refer the reader to \cite[\S 4]{Mil74} for the details. 

Let $\xi: E\rightarrow B$ denote a vector bundle over a CW complex $B$ and  $w_i(\xi)\in H^i(B,\mathbb{Z}_2)$ denote the $i^{th}$ Stiefel-Whitney class of $\xi$. If $\xi$ is a trivial bundle, then $w_i$ vanishes for all $i$. The first Stiefel-Whitney class detects the orientability of $\xi$. That is,  $w_1(\xi)=0$ if and only if $\xi$ is orientable. Suppose that $\xi$ is an oriented vector bundle of dimension $n$, then the $n^{th}$ Stiefel-Whitney class $w_n\in H^n(B,\mathbb{Z}_2)$ is the mod $2$ reduction of the Euler class $e(\xi)$ in $H^n(B,\mathbb{Z})$. When $n=2$, if we use the notation introduced in \S \ref{subsec: define euler class}, this means that  $w_2(\xi)=[\bar{c}_\sigma]$ with $\bar{c}_\sigma(e_\alpha)=c_\sigma(e_\alpha) \mbox{ (mod $2$)}$ for each $2$-cell $e_\alpha$ of $B$. Lastly, if $\xi\cong \xi_1\oplus \xi_2$ is the Whitney sum of two vector bundles  $\xi_1: E_1\rightarrow B$ and $\xi_2: E_2\rightarrow B$, then we have  $w_i(\xi)=\sum_{k=0}^i w_k(\xi_1)w_{i-k}(\xi_2)$.  

It is a well-known fact that the tangent bundle $TM$ of a closed orientable $3$-manifold $M$ is trivial. To prove Proposition \ref{prop:Euler class is even}, we only need $w_2(TM)=0$, which can be easily verified using Wu's formula \cite[\S 11]{Mil74}. Also see \cite[Theorem 4.2.1]{Geiges2008} 

\begin{proof}[Proof of Proposition \ref{prop:Euler class is even}]
Since both $\xi$ and $TM$ are orientable, the normal line bundle $\eta=TM/\xi$ is orientable and hence trivial.  From the decomposition $TM\cong \eta\oplus \xi$, we have 
$$w_2(\xi)=w_2(TM)-w_1(\xi)w_1(\eta)-w_2(\eta)=w_2(TM).$$ As $w_2(TM)=0$, so is $w_2(\xi)$.

Consider the Bockstein sequence associated to the short exact sequence $0 \rightarrow \mathbb{Z}\xrightarrow{\cdot 2} \mathbb{Z} \rightarrow \mathbb{Z}_2\rightarrow 0$ as follows:
\begin{equation}
\rightarrow  H^2(M,\mathbb{Z})\xrightarrow{\cdot 2} H^2(M,\mathbb{Z})\xrightarrow{{\rm mod} \, 2} H^2(M,\mathbb{Z}_2)\rightarrow H^3(M,\mathbb{Z})\rightarrow.  
\label{equ:mod 2}
\end{equation}
Since the Euler class $e(\xi)$ is sent to $w_2(\xi)=0$ under $H^2(M,\mathbb{Z})\xrightarrow{{\rm mod} \, 2} H^2(M,\mathbb{Z}_2)$, we have  $e(\xi)$ is in the image of  $ H^2(M,\mathbb{Z})\xrightarrow{\cdot 2} H^2(M,\mathbb{Z})$ as claimed.

Suppose that $H^2(M,\mathbb{Z})\cong \oplus \mathbb{Z}_2$ is a direct sum of $\mathbb{Z}_2$. Then the homomorphism $H^2(M,\mathbb{Z})\xrightarrow{\cdot 2} H^2(M,\mathbb{Z})$ is the zero map. Hence, the Euler class $e(\xi) \in {\rm Im}(H^2(M,\mathbb{Z})\xrightarrow{\cdot 2} H^2(M,\mathbb{Z}))$ must be zero. 
\end{proof}

\section{Necessary and sufficient conditions for the Euler class to vanish} 
\label{sec: euler class and dehn filling}
In this section, we prove Theorem \ref{thm: e=0 necessary sufficient condition}. It gives necessary and sufficient conditions for the Euler class of a foliation on a Dehn filling of a $\mathbb{Q}$-homology solid torus $X$ to vanish. Theorem \ref{thm: euler class zero knots in zhs} is the special case of Theorem \ref{thm: e=0 necessary sufficient condition} when $X$ is a $\mathbb{Z}$-homology solid torus. 

\subsection{Notation and conventions}
\label{subsec: notation and orientation}
Let $X$ be a $\mathbb{Q}$-homology solid torus. That is, $X$ is a compact orientable $3$-manifold with $H_\ast(X;\mathbb{Q})\cong H_\ast(D^2\times S^1; \mathbb{Q})$. It follows that $\partial X$ consists of a single torus. Let $F$ be an oriented connected incompressible, boundary incompressible surface properly embedded in $X$ such that $[F]$ is a generator of $H_2(X,\partial X)\cong \mathbb{Z}$. We may assume that there are $k\geq 1$ connected components of $\partial F$ on $\partial X$ oriented coherently by the induced orientation from $F$. The {\it rational longitude} $\lambda$ of $X$, or {\it longitude} for short,  is a connected component of $\partial F$.  Let $\iota: \partial X\rightarrow X$ be the inclusion map. Then $\lambda$ represents the unique primitive integral class up to sign that is in the kernel of  $\iota_*: H_1(\partial X,\mathbb{Q})\rightarrow H_1(X, \mathbb{Q})$. By definition, $k$ is the order of $\iota_*([\lambda])$ in $H_1(X,\mathbb{Z})$. We refer to $k$ as {\it the order of $\lambda$}.

A {\it meridian} $\mu$ is an oriented simple closed curve on $\partial X$ that intersects $\lambda$ once. Hence, $\mu$ and $\lambda$ form a basis of $H_1(\partial X)$. We orient $\mu$ so that the algebraic intersection number of $F$ and $\mu$ equals $k$. In particular, $\mu$ intersects both $F$ and $\lambda$ positively.  In Figure \ref{fig:unknot frame}, we use the exterior of the unknot $U$ in $S^3$ with the standard orientation (the right hand rule) to illustrate our choice of the orientations on $\lambda$, $\mu$ and $F$.

\begin{figure}[ht]
\begin{center}
\begin{tikzpicture}[scale=0.7]
 %\draw [help lines] (0,0) grid (12, 6);
 \draw (6,3) circle [radius=1.5];
 \node [right] at (6.9, 4) {\footnotesize $U$};
 \draw [fill, light-gray, thick](6,3) circle [radius=1.2];
  \begin{scope}[decoration={
    markings,
    mark=at position 0.18 with {\arrow{>}}}] 
 \draw [blue, thick, postaction={decorate}] (7.2,3) to [out=90, in=0] (6,4.2) to [out=180, in=90] (4.8,3) to [out=270, in=180] (6,1.8) to [out=0, in=270] (7.2,3); 
\end{scope}
 \draw [thick](6,3) circle [radius=1.8];
 \begin{scope}[decoration={
    markings,
    mark=at position 0.5 with {\arrow{>}}}] 
 \draw [red, thick, postaction={decorate}] (7.19,3) to [out=80, in=180] (7.5, 3.3) to [out=0, in=100] (7.81,3);
\end{scope}
\draw [red, dashed] (7.19,3) to [out=-80, in=180] (7.5, 2.8) to [out=0, in=-100] (7.81,3);
\node [red]at (7.6, 3) {\footnotesize $\mu$};
\node [blue] at (6, 4) {\footnotesize $\lambda$};
\node at (6,3) {\footnotesize $F\quad +$ };
\end{tikzpicture}
\end{center}
\caption{}
\label{fig:unknot frame}
\end{figure}

If $X$ is the exterior of a null-homologous knot in a $\mathbb{Q}$-homology sphere, then the boundary of a meridional disk of the tubular neighborhood of the knot is a meridional curve, called {\it the knot meridian}.

A slope $\alpha$ on $\partial X$ is an isotopy class of essential simple closed curves on $\partial X$. Given $\mu$ and $\lambda$ as above, we write $\alpha=p/q$ in $\mathbb{Q} \cup \{\frac{1}{0}\}$, if  $[\alpha]=\pm (p[\mu]+q[\lambda])$ in $H_1(\partial X)$. So $\lambda =0$ and $\mu = \frac{1}{0}$. If we choose a different meridian $\mu'$ of $X$, then since $[\mu'] = [\mu] + n[\lambda]$ in $H_1(\partial X)$ for some $n\in \mathbb{Z}$, we have $\alpha$, which equals $p/q$ with respect to $\mu$, now is equal to $p/(q-np)$ with respect to $\mu'$. Given an oriented simple closed curve $c$ on $\partial X$, to simplify the notation, we often use $c$ to denote both the slope and the homology class in $H_1(\partial X)$ it represents.

We use $X(\alpha)$ to denote the closed $3$-manifold obtained by the $\alpha$-Dehn filling on $X$. That is, $X(\alpha)=X\cup_f N$ is obtained by attaching a solid torus $N$ to $X$ along $\partial X$ so that the boundary of the meridional disk $D$ of $N$ is identified with $\alpha$ under the gluing map $f : \partial N\rightarrow \partial X$. We orient the disk $D$ so that $f(\partial D)=p\mu+q\lambda$ with $p\geq0$.

\subsection{Cohomology classes in $H^2(X(\alpha))$}
In $X(\alpha)$, we identify $\partial N$ with $\partial X$ by the map $f$. We assume that $\alpha\neq \lambda$, so $X(\alpha)$ is a $\mathbb{Q}$-homology sphere and we have the following diagram:

\begin{center}
\begin{tikzpicture} [scale=0.8]
 %\draw [help lines] (0,0) grid (12, 6);
\node at (6, 4) {\small $ 0 \longrightarrow H^1(\partial X)\xrightarrow{\;\;\; \delta \;\;} H^2(X(\alpha),\partial X)\xrightarrow{\;\;\; i^* \;\;} H^2(X(\alpha))\longrightarrow 0 $};
\draw [->] (5.8, 2.7) -- (5.8, 3.5); 
\draw [->] (5.8, 1.3) -- (5.8, 1.9);
\draw [->] (5.8, 4.3) -- (5.8, 4.9);
\node at (5.8, 5.2) {\small $0$};
\node at (5.8, 1) {\small $0$};
\node at (5.8, 2.3) {\small $H^2(X,\partial X)\oplus H^2(N,\partial N)$};
\node [right] at (5.8, 3.1) {\small $\cong$};
\end{tikzpicture}
\end{center}
where the vertical isomorphism is from the decomposition 
\begin{displaymath}
(X(\alpha), \partial X) = (X,\partial X) \cup (N, \partial N = f^{-1}(\partial X)).
\end{displaymath}

More precisely, given a cohomology class $c\in H^2(X(\alpha), \partial X)$, the corresponding class in $H^2(X,\partial X)\oplus H^2(N,\partial N)$ is given by $(c|_X, c|_N)$, where $c|_X$ and $c|_N$ are restrictions of $c$ to the subspaces $X$ and $N$ respectively. On the other hand, given relative classes $c_X$ and $c_N$ in $H^2(X,\partial X)$ and $H^2(N, \partial N)$, since both of them vanish on the boundaries, we may extend them to the closed manifold $X(\alpha)$ by the zero map.  We continue to use  $c_X$ and $c_N$ to denote the resulting classes in $H^2(X(\alpha), \partial X)$. Then under the isomorphism 
$H^2(X,\partial X)\oplus H^2(N,\partial N) \cong H^2(X(\alpha), \partial X)$, the pair $(c_X, c_N)$ is mapped to  the sum $c_X +c_N$.

By identifying $H^2(X(\alpha), \partial X)$ and $H^2(X,\partial X)\oplus H^2(N,\partial N)$ as above, we obtain the short exact sequence 
\begin{equation}
  0 \longrightarrow H^1(\partial X)\xrightarrow{\;\;\; \delta \;\;} H^2(X,\partial X)\oplus H^2(N,\partial N)\xrightarrow{\;\;\; i^* \;\;} H^2(X(\alpha))\longrightarrow 0. 
  \label{equ: exact sequence}
\end{equation}
where $i^*( c_X ,c_N) = c_X + c_N$ and $\delta \beta= (\delta \beta, \delta\circ f^*\beta)$ for any $(c_X, c_N)\in H^2(X,\partial X)\oplus H^2(N,\partial N)$ and $\beta\in H^1(\partial X)$. Notice that the short exact sequence (\ref{equ: exact sequence}) is the dual of the Mayer-Vietoris sequence of homology groups given by the decomposition $X(\alpha) = X\cup N$.
 
The lemma below contains the core of the computation in this paper. Later in \S\ref{subsec: euler classes taut foliation}, the cohomology class $c$ in the lemma will be the Euler class $e(\widehat{\mathcal{F}})$, and $c'|_X$ and $c'|_N$ will be two relative Euler classes of $\mathcal{F} = \widehat{\mathcal{F}}|_X$ and $\mathcal{D} = \widehat{\mathcal{F}}|_N$ respectively (see Theorem \ref{thm: e=0 necessary sufficient condition}).

\begin{lemma} 
Fix a meridian $\mu$ and let $\alpha=p/q$ with $p>0$. Assume that $c\in H^2(X(\alpha))$ and $c'$ is a relative class in $H^2(X(\alpha),\partial X)$ satisfying $i^*(c')=c$. Let $c'|_X$ and $c'|_N$ denote the restrictions of $c'$ to $H^2(X,\partial X)$ and $H^2(N,\partial N)$ respectively. Suppose that $c = 0$ in $H^2(X(\alpha))$. Then the following two statements hold:
\begin{enumerate}
    \item The restriction $c|_X = 0$ in $H^2(X)$;
    \item Suppose that  $a = c'|_X([F])$, $b= c'|_N([D])$ and  $k$ is the order of the rational longitude $\lambda$. Then $a$ is divisible by $k$ and $(aq)/k  \equiv b \mbox{ \rm{(mod $p$)}}$. 
\end{enumerate}
In addition, if $k=1$, i.e., $\lambda$ is null-homologous in $H_1(X)$, then the above two conditions are sufficient to conclude that $c = 0$.
 \label{lem:cohomology class zero}
\end{lemma}

\begin{remark}
The statement of Lemma \ref{lem:cohomology class zero} doesn't depend on the choice of meridian $\mu$. However, it does depend on the choice of orientations. A different orientation convention than the one we described in \S \ref{subsec: notation and orientation} may result in $(aq)/k\equiv -b \mbox{ (mod $p$)}$ in Condition (2).  
\label{rem: dependence on the orientation lemma 3.3}
\end{remark}

\begin{proof}
The necessity of Condition (1) is obvious. Since the restriction $c|_X$ is the image of $c$ under the inclusion induced map  $H^2(X(\alpha)) \rightarrow H^2(X)$, it's zero if $c$ is. 

Next we derive Condition (2) by assuming $c = 0$. Since $c= i^*(c')=i^*(c'|_X, c'|_N)=0$, there exists $\beta\in H^1(\partial X)$ such that $(c'|_X, c'|_N) = (\delta\beta, \delta\circ f^*\beta)$. That is, 
\begin{equation}
     c'|_X = \delta\beta\, \text{ and } \,c'|_N= \delta\circ f^* \beta. 
     \label{equ: 0}
\end{equation}
Hence,  
\begin{displaymath}
   b = c'|_N([D]) = \delta\circ f^*\beta ([D]) = \beta (f(\partial D)) = q\beta(\lambda) + p\beta(\mu),
\end{displaymath}
whereas in the first term 
\begin{equation}
\label{equ: a divisible by k}
  \beta(\lambda) = \delta\beta([F])/k = c'|_X([F])/k = a/k. 
\end{equation}
Therefore,  $a = c'|_X([F])=k\beta(\lambda)$ is divisible by $k$ and $(aq)/k  \equiv b \mbox{ (mod $p$)}$ as claimed. 

\medskip

Now assuming that (1) and (2) hold and $k=1$, we show that $c=0$ in $H^2(X(\alpha))$. 
Applying $H^2$ to the following commutative diagram: 
\begin{center}
\begin{tikzpicture}[scale=0.9]
%\draw [help lines] (0,0) grid (12, 6);
 \node at (7,4) {\small $(X,\partial X)$};
 \draw [thick, ->] (8, 4)--(9,4);
 \node at (10.3,4) {\small $(X(\alpha), \partial X)$};
 \draw [thick, <-] (7.4, 3.6) -- (7.4, 2.6);
 \node at (7.1, 3.1) {$i$};
 \node at (9.8, 2.2) {\small $X(\alpha)$};
 \draw [thick, ->] (8, 2.2)--(9,2.2);
 \node at (7.4, 2.2) {\small $X$};
 \draw [thick, <-] (9.7, 3.6) -- (9.7, 2.6);
  \node at (10, 3.1) {$i$};
  %\node at (8.5, 2.4) {\small res. };
  %\node at (8.5, 4.2) {\small res. };
\end{tikzpicture}
\end{center}
we have $i^*(c'|_X) = c|_X$ in $H^2(X)$. Since $c|_X = 0$ by assumption, from the exact sequence
\begin{displaymath}
 H^1(X) \xrightarrow{\,\iota^*} H^1(\partial X) \xrightarrow{\delta} H^2(X, \partial X) \xrightarrow{\,i^*} H^2(X), 
\end{displaymath}
there exists $\beta_0\in H^1(\partial X)$ such that $\delta \beta_0 = c'|_X.$ We fix $\beta_0$. Given any $\beta' \in H^1(\partial X)$ with $\delta \beta' = c'|_X$, the difference $\beta_0 - \beta'$ is in the image of $\iota^*$, which is generated by $k\mu^* = \mu^*$ (since $k=1$). Here $\mu^*$ in $H^1(\partial X)$ is the dual of $\mu$. Let $\beta = \beta_0 + n\mu^*$ for some $n\in\mathbb{Z}$.  To show that $c = 0$, it remains to show that there exists $n\in \mathbb{Z}$ such that  $\delta\circ f^*\beta =  c'|_N$ in $H^2(N, \partial N)$. 

Since $H^2(N,\partial N)\cong {\rm Hom}(H_2(N,\partial N),\mathbb{Z})$, this is equivalent to $\delta\circ f^*\beta([D]) =  c'|_N([D])$. By Condition (2),  $aq - b$ is divisible by $p$. Let $n = -\beta_0(\mu) - \frac{aq - b}{p}$ and $\beta = \beta_0 + n\mu^*$. Then 
 \begin{align*}    
   \delta\circ f^*\beta ([D]) & = \beta (f(\partial D)) \\ \nonumber
     & = (\beta_0 + n\mu^*)(p\mu+q\lambda) \\ 	
     & =  p\beta_0(\mu) + q\beta_0(\lambda) + np \\
     & = p\beta_0(\mu) + q\beta_0(\lambda)  + \left(-\beta_0(\mu) - \frac{aq - b}{p}\right)p \\
     & = q\beta_0(\lambda) - aq + b \\
    & = q\delta\beta_0([F]) - aq +b \quad \text{(since $k=1$)}\\
     &= qc'|_X([F]) - aq+b \\
     & = b =c'|_N([D]).
\end{align*}
Hence 
we have $\delta\beta = (\delta \beta, \delta\circ f^*\beta)=(c'|_X, c'|_N)= c'$. Therefore,  $c = i^*(c') = 0$. 
\end{proof}

\begin{remark}
\label{rem: iff condition k>1}
In general, the argument in the proof of Lemma \ref{lem:cohomology class zero} shows that a necessary and sufficient condition for $c=0$ is that both  $c|_X = 0$  and  $(aq)/k + \beta(\mu)p \equiv b \mbox{ \rm{(mod $kp$)}}$ hold. Here $\beta\in H^1(\partial X)$ satisfies $\delta\beta = c'|_X$ whose existence is guaranteed by the condition $c|_X = 0$.  As we have seen in the proof of Lemma \ref{lem:cohomology class zero},  the value $(aq)/k + \beta(\mu)p \mbox{ \rm{(mod $kp$)}}$  doesn't depend on the choice of $\beta$.

\end{remark}

\subsection{The Euler class of foliations on Dehn filled manifolds.}
\label{subsec: euler classes taut foliation}
In this section, we apply Lemma \ref{lem:cohomology class zero} to computing the Euler class of foliations on $X(p/q)$, $p>0$.  
Let $\widehat{\mathcal{F}}$ denote an oriented co-dimension-one foliation on $X(\alpha)$ such that $\mathcal{D} = \widehat{\mathcal{F}}|_N$ is the foliation by meridional disks. Consequently, $\mathcal{F} = \widehat{\mathcal{F}}|_X$ intersects  $\partial X$ transversely in parallel simple closed curves of slope $\alpha$. 

Let $T\widehat{\mathcal{F}}$, $T\mathcal{F}$ and $T\mathcal{D}$ denote the oriented tangent plane fields of $\widehat{\mathcal{F}}$, $\mathcal{F}$ and $\mathcal{D}$ respectively. We use $e(\widehat{\mathcal{F}})$ to denote the Euler class of $T\widehat{\mathcal{F}}$  in $H^2(X(p/q))$. Let $\sigma$ be a nowhere vanishing section of $T\mathcal{F}|_{\partial X}=T\widehat{\mathcal{F}}|_{\partial X}$ that is transverse to $\partial X$ pointing outwards. Hence, $\sigma$ is also a nowhere vanishing section of $T\mathcal{D}$ along $\partial N$ pointing into $N$.

Following \S \ref{subsec: define euler class}, we obtain three relative Euler classes. They are $e_\sigma(\widehat{\mathcal{F}})$ in $H^2(X(\alpha), \partial X)$,  $e_\sigma(\mathcal{F})$ in $H^2(X,\partial X)$ and $e_\sigma(\mathcal{D})$ in $H^2(N, \partial N)$.  By definition, 
$e_\sigma(\widehat{\mathcal{F}})|_X=e_\sigma(\mathcal{F})$ and $e_\sigma(\widehat{\mathcal{F}})|_N=e_\sigma(\mathcal{D})$, where $e_\sigma(\widehat{\mathcal{F}})|_X$ and $e_\sigma(\widehat{\mathcal{F}})|_N$ are restrictions of $e_\sigma(\widehat{\mathcal{F}})$ to  $X$ and $N$ respectively. Moreover,  $i^*(e_\sigma(\widehat{\mathcal{F}}))=e(\widehat{\mathcal{F}})$ and $i^*(e_\sigma(\mathcal{F})) = e(\mathcal{F})$, where $i^*$'s are induced by inclusions $(X(p/q),\emptyset) \rightarrow (X(p/q), \partial X)$ and $(X,\emptyset) \rightarrow (X, \partial X)$ respectively.

Theorem \ref{thm: e=0 necessary sufficient condition} below follows immediately from Lemma \ref{lem:cohomology class zero}. Euler classes $e(\widehat{\mathcal{F}})$,  $e(\mathcal{F})$, $e_\sigma(\widehat{\mathcal{F}})$, $e_\sigma(\mathcal{F})$ and $e_\sigma(\mathcal{D})$ correspond to cohomology classes $c$, $c|_X$, $c'$, $c'|_X$ and $c'|_N$ in the lemma.

\begin{theorem}
Let $X$ be a $\mathbb{Q}$-homology solid torus and $F$ be a properly embedded surface in $X$ representing a generator of  $H_2(X,\partial X)$. Suppose that $\widehat{\mathcal{F}}$ is an oriented co-dimension one foliation on $X(p/q)$, $p>0$ whose restriction  to the filling solid torus $N$ is  the foliation by meridian disks, and the orientation of the leaves of $\widehat{\mathcal{F}}$ agree with the given orientation of the meridian disks of $N$. Let $\mathcal{F} = \widehat{\mathcal{F}}|_X$ and $\sigma$ denote a nowhere vanishing outward pointing section of $T\mathcal{F}$ along $\partial X$.  Suppose that the Euler class $e(\widehat{\mathcal{F}}) = 0$ in $H^2(X(p/q))$. Then the following two statements hold: 
\begin{enumerate}
    \item The Euler class $e(\mathcal{F}) = 0$ in $H^2(X)$;
    \item Let $a = e_\sigma(\mathcal{F})([F])$ and $k$ denote the order of the rational longitude $\lambda$. Then $a$ is divisible by $k$ and $(aq)/k  \equiv 1 \mbox{ \rm{(mod $p$)}}$.
\end{enumerate}
In addition, if $k=1$, i.e., $\lambda$ is null-homologous in $H_1(X)$, then the above two conditions are sufficient to conclude that $e(\widehat{\mathcal{F}}) = 0$.
\label{thm: e=0 necessary sufficient condition}
\end{theorem} 

\begin{remark}
\label{rem: e(F) = 0 iff condition k>1}
A necessary and sufficient condition for $e(\widehat{\mathcal{F}}) =0$ when $k> 1$ in general can be deduced immediately from the condition given in Remark \ref{rem: iff condition k>1} with $b$ being replaced by $1$. 
\end{remark}

Intuitively, Condition (1) in Theorem \ref{thm: e=0 necessary sufficient condition} says that the plane field $T\widehat{\mathcal{F}}$ is trivial over $X$, which is obviously necessary for it to be trivial over the entire manifold $X(\alpha)$. Also notice that since $H^2(N)=0$, we have $T\widehat{\mathcal{F}}|_N$ is always trivial. Therefore,  there exist nowhere vanishing sections of $T\widehat{\mathcal{F}}$ over both $X$ and $N$. Condition (2) spells out the numerical equation that must be satisfied for the existence of sections over $X$ and $N$ that match along $\partial X$, so that  together they define a  global section of $T\widehat{\mathcal{F}}$ over $X(\alpha)$.

\begin{proof}[Proof of Theorem \ref{thm: e=0 necessary sufficient condition}]
The theorem is Lemma \ref{lem:cohomology class zero} with an additional claim that the value of $b$  must equal $1$.

By assumption $\mathcal{D} := \widehat{\mathcal{F}}|_N$ is the foliation of $N$ by meridional disks. Let $D$ denote a meridian disk of $N$. Since the orientation of $D$ inherited from $\mathcal{D}$ agrees with the given orientation on $D$ (see \S \ref{subsec: notation and orientation}),  by the Poincare-Hopf Theorem, we have $b=e_\sigma(\mathcal{D})([D])=\chi(D)= 1.$ 
\end{proof}

\begin{proof}[Proof of Theorem  \ref{thm: euler class zero knots in zhs}]
When $X$ is a $\mathbb{Z}$-homology solid torus, $H^2(X) = 0$ and $k = 1$.  Therefore, Theorem  \ref{thm: euler class zero knots in zhs} follows from Theorem \ref{thm: e=0 necessary sufficient condition}.
\end{proof}

\section{Dehn fillings on fibered manifolds and left-orders}
\label{sec: examples}
In this section, we compute the Euler classes of taut foliations constructed in \cite{RobertsSurfacebundle2} (Theorem \ref{thm: euler class Rachel's foliation}) and \cite{krishna18} (Theorem \ref{thm: foliation 3 braids}). These are taut foliations on closed $3$-manifolds obtained from Dehn fillings on fibered manifolds.  

\subsection{Canonical meridian, degeneracy slope and the fractional Dehn twist coefficient}
\label{subsec: fibered knot}
Let $X$ be the exterior of a fibered hyperbolic knot with fiber $F$ in a $\mathbb{Q}$-homology sphere. So $X=F\times [0,1]/(x,1)\sim_h (h(x), 0)$ for an orientation-preserving homeomorphism $h: F\rightarrow F$ with $h|_{\partial F} = {\rm id}_{\partial F}$, where $h$ is called the monodromy of the knot. Since the interior of $X$ is hyperbolic, $h$ is freely isotopic to a pseudo-Anosov homeomorphism $\varphi$ \cite{Thurston98}. Note that $\varphi|_{\partial X}$ is never the identity. 

The flow lines of the suspension flow of $\varphi$ on $\partial X$ are path-connected components of the quotient $\sqcup_{p\in \partial F} (p\times [0,1]) /\sim_\varphi$. The {\it degeneracy slope} of $X$, denoted by $\gamma$ is a closed flow line of the suspension flow of $\varphi$, the existence of which is guaranteed by the properties of pseudo-Anosov maps (see \cite[Part 3]{FM} for instance). The concept of degeneracy slope was first considered in \cite{GO89} in the study of essential laminations in Dehn fillings of $X$. 

Instead of the knot meridian $\mu=\ast\times [0,1]/\sim_h$, when considering fibered knots in a general $\mathbb{Q}$-homology sphere, it is often natural to use the so-called {\it canonical meridian} \cite{{RobertsSurfacebundle2}}, which we denote by $\mu_0$. 

To obtain $\mu_0$, we first follow $\gamma$ starting at a point in $\gamma\cap \lambda$ until we reach a point in $\gamma\cap \lambda$ again. In general, this is not the same intersection point as the one we started with. To form a loop, which represents the slope $\mu_0$, we continue to follow one of the subarcs of the longitude $\lambda$ back to the initial point. There are two subarcs of $\lambda $ that can lead us back to the initial intersection point, and we choose the one that intersects $\gamma$ minimally.  If both subarcs have the same intersection number with $\gamma$, we follow the convention in \cite[\S 3]{RobertsSurfacebundle2} and use the subarc of $\lambda$ so that $\gamma$ is a positive slope with respect to $\mu_0$. 

\begin{remark}[The trace of the monodromy of genus one knots]
\label{rem: genus one trace}
Let $h$ denote the monodromy of a fibered knot $K$ and $h_*: H_1(F)\rightarrow H_1(F)$ be its induced homology map. When the genus  $g(F) = 1$, we have $\gamma = \mu_0$ if and only if the trace ${\rm tr}(h_*) > 2$. This is because when $g(F)=1$, the stable foliation $\mathcal{F}_s$ of the (pseudo-)Anosov representative $\varphi$ of $h$ has two singular points on the boundary. The orientation of $\mathcal{F}_s$ is preserved by $\varphi$ when $\rm{tr}(h_*)>2$, so $\varphi$ must fix both singular points of $\mathcal{F}_s$ on $\partial F$. Therefore, the degeneracy slope $\gamma$ intersects $\lambda$ once and hence equals  the canonical meridian $\mu_0$.  
\end{remark}

There is another notion closely related to the degeneracy slope, called the {\it fractional Dehn twist coefficient} $c(h)$ of $h$. It was first introduced in \cite{HKMII,HKMI} to study the tightness of the contact structure supported by the open book $(F, h)$.  If we write the degeneracy slope $\gamma = p\mu + q\lambda$, where $\mu$ is the knot meridian, then the fractional Dehn twist coefficient $h$, denoted by $c(h)$, equals $q/p$, the reciprocal of the degeneracy slope \cite[\S 2]{KR13}.

\begin{remark}[The degeneracy slope and the FDTC]
\label{rem: FDTC and degeneracy slope}
\begin{enumerate}
\item  When $c(h)$ is an integer, the degeneracy slope $\gamma$ intersects $\lambda$ once. Therefore, by definition $\mu_0 = \gamma $. When $X$ is the exterior of a fibered knot in $S^3$, it is known that $|c(h)|<1$. Examples of fibered knots in $S^3$ with $c(h) = 0$ include hyperbolic fibered two-bridge knots \cite[Theorem 8]{GK90}.

\item When $c(h) > 0$,  $h$ is called {\it right-veering}. Fibered knots in $S^3$ with right-veering monodromies include the positive L-space knots and more generally fibered strongly quasipositive knots \cite{Ni07, Ghi08, HKMI, Hed}. For these knots, the degeneracy slope is positive with respect to the knot meridian and hence the canonical meridian $\mu_0$ is the same as the knot meridian  \cite[Proposition 7.4]{RobertsSurfacebundle2}. 
\end{enumerate}
 \end{remark}

\subsection{Branched surfaces}
\label{subsec: euler class branched surfaces}
\begin{figure}[ht]
\centering
\begin{tikzpicture}[scale=0.65]
%\draw [help lines] (0,0) grid (18,5);
\draw (1,1) -- (4,1) -- (5,3) -- (2,3) -- (1,1);
% the middle 
\draw (7,1) -- (10,1) -- (11,3) -- (8,3) -- (7,1);
\draw (8.5,1) to [out=0, in=235] (9.5, 1.5) to  (10.5,3.5);
\draw (8.5,1) -- (9.5, 3);
\draw (9.5, 3) to [out=0, in=235] (10.5,3.5);
\draw [-{>[scale=2.5, length=1, width=1]}] (9.2, 2) -- (8.7, 2);
% the third one
\draw (13,1) -- (16,1) -- (17,3) -- (14,3) -- (13,1);
\draw (14.5,1) to [out=0, in=235] (15.5, 1.5) to  (16.5,3.5);
\draw (14.5,1) -- (15.5, 3);
\draw (15.5, 3) to [out=0, in=235] (16.5,3.5);
\draw [-{>[scale=2.5, length=1, width=1]}] (14.9, 1.5) -- (14.5, 1.5);
\draw [-{>[scale=2.5, length=1, width=1]}] (14.3, 1.85) -- (14.3, 2.2);
%% the bottom piece
\draw [dashed] (13.5,2) -- (16.5,2);
\draw [dashed] (13.4,1.15) -- (16.05,1.15);
\draw (16.05,1.15) -- (16.4,1.15);
\draw (16.5,2) to [out = 250, in=85] (16.4,1.15);
\draw [dashed] (13.5,2) to [out = 250, in=85] (13.3,1.15);
\end{tikzpicture}
\caption{Local models of a branched surface}
\label{fig:branched surface}
\end{figure}
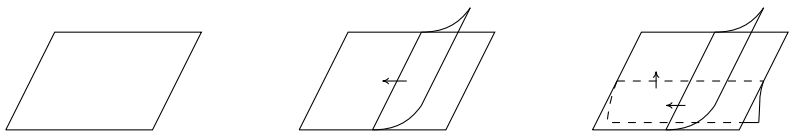

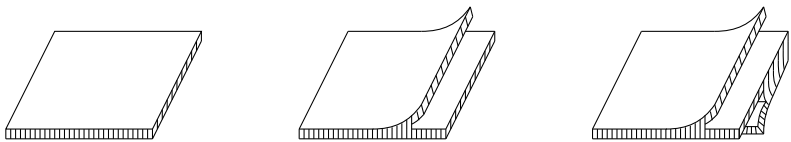
\begin{figure}[ht]
\centering
\begin{tikzpicture}[scale=0.65]
%\draw [help lines] (0,0) grid (18,5);
\draw (1,1) -- (4,1) -- (5,3) -- (2,3) -- (1,1);
\draw (1, 0.8) -- (4, 0.8);
\draw (4, 0.8) -- (5,2.8);
\foreach \x in {1, 1.1, 1.2, ..., 4}
	\draw (\x, 0.8) -- (\x, 1);
\foreach \n in {1,2,...,10}
	\draw (4+0.1*\n, 0.8+0.2*\n) -- (4+0.1*\n, 1+0.2*\n);
% the middle 
%\draw (7,1) -- (10,1) -- (11,3) -- (8,3) -- (7,1);
\draw (10,1) -- (11,3);
\draw (8,3) -- (7,1);
\draw (7,1) -- (8.5,1);
\draw (8,3) -- (9.5,3);
\draw (9.3,1) -- (10,1);
\draw (8.5,1) to [out=0, in=235] (9.5, 1.5) to  (10.5,3.5);
\draw (9.5, 3) to [out=0, in=235] (10.5,3.5);
%% middle
\draw (7, 0.8) -- (10, 0.8);
\draw (10, 0.8) -- (11,2.8);
\foreach \x in {1, 1.1, 1.2, ..., 4}
	\draw (6+\x, 0.8) -- (6+\x, 1);
\foreach \n in {1,2,...,10}
	\draw (10+0.1*\n, 0.8+0.2*\n) -- (10+0.1*\n, 1+0.2*\n);
%% front part
\foreach \n in {0,1,2,...,7,8}
	\draw (8.5+0.1*\n, 1) -- (8.5+0.1*\n, 1+0.00055*\n*\n*\n);
\draw (9.3, 1.1) to [out=25, in=225] (9.55, 1.3);
\draw (9.45, 1.18) -- (9.41, 1.36);
%% top layer
\draw (10.4, 3) -- (11,3);
\draw (9.55, 1.3) --  (10.55,3.3);
\foreach \n in {0,1,2,...,10}
	\draw (9.55 + 0.1*\n , 1.3+ 0.2*\n) -- (9.5+0.1*\n , 1.5+ 0.2*\n);
% the third one
\draw (16,1) -- (17,3);
\draw (14,3) -- (13,1);
\draw (13,1) -- (14.5,1);
\draw (14,3) -- (15.5,3);
\draw (15.3,1) -- (16,1);
\draw (14.5,1) to [out=0, in=235] (15.5, 1.5) to  (16.5,3.5);
\draw (15.5, 3) to [out=0, in=235] (16.5,3.5);
%% front part
\foreach \n in {0,1,2,...,7,8}
	\draw (14.5+0.1*\n, 1) -- (14.5+0.1*\n, 1+0.00055*\n*\n*\n);
\draw (15.3, 1.1) to [out=25, in=225] (15.55, 1.3);
\draw (15.45, 1.18) -- (15.41, 1.36);
%% top layer
\draw (16.4, 3) -- (17,3);
\draw (15.55, 1.3) --  (16.55,3.3);
\foreach \n in {0,1,2,...,10}
	\draw (15.55 + 0.1*\n , 1.3+ 0.2*\n) -- (15.5+0.1*\n , 1.5+ 0.2*\n);
%% middle layer
\draw (13, 0.8) -- (16, 0.8);
\draw (16, 0.8) -- (16.5, 1.8);
\foreach \x in {1, 1.1, 1.2, ..., 4}
	\draw (12+\x, 0.8) -- (12+\x, 1);
\draw (17,2.4) -- (16.7, 1.7) to [out=245, in=90] (16.5,0.9);
\foreach \n in {0, 1, 2, ...,10}
	\draw (16+ 0.1*\n, 1+0.2*\n) to (16+ 0.1*\n, 0.8+0.2*\n);
\foreach \n in {0, 1, 2, ...,5}
	\draw (16.5+ 0.1*\n, 1.8+0.2*\n) to [out = 270, in = 135-9*\n] (16.6+0.08*\n, 1.5+0.18*\n);
%% right bottom part
\draw (16.38, 1.56) -- (16.55, 1.56);
\draw (16.05,0.9) -- (16.5, 0.9);
\draw (16.55, 1.56) to [out= 235, in =85] (16.38,1.05) -- (16.12, 1.05);
%%% the fibers manually front
\draw (16.38, 1.05) to (16.5,0.9);
\draw (16.25, 1.05) -- (16.25, 0.9);
\draw (16.17, 1.05) -- (16.15, 0.9);
\draw (16.33, 1.05) -- (16.35, 0.9);
%%% the fibers manually side
\draw (16.39, 1.16) -- (16.51, 1.06);
\draw (16.41, 1.28) -- (16.52, 1.19);
\draw (16.45, 1.37) -- (16.55, 1.29);
\draw (16.50, 1.49) -- (16.58, 1.4);
\end{tikzpicture}
\caption{Local models of the normal neighborhood of a branched surface}
\label{fig:normal neighborhood}
\end{figure}

A branched surface $B$ is a topological space locally modeled on one of the pictures in  Figure \ref{fig:branched surface}. The set of points at which $B$ is not diffeomorphic to $\mathbb{R}^2$ is called the branch locus. Note that at any point of $B$ there is a well-defined tangent plane. Given a properly embedded branched surface $B$ in a $3$-manifold $M$, let $N(B)$ denote an $I$-fibered normal neighborhood of $B$ in $M$, whose local models are depicted in Figure \ref{fig:normal neighborhood}.  

\begin{example}
\label{ex: branched surface fd}
Let $X$ be the exterior of a fibered knot with oriented fiber $F$ and monodromy $h$. Given a properly embedded arc $\beta$ in $F$, one can construct a branched surface $B = \langle F;\mathbb{D}\rangle$ in $X$ by attaching a disk $\mathbb{D} = [0,1]\times [0,1]$ to $F$ with the bottom side $[0,1]\times 0$ of $\mathbb{D}$ attached to the positive side of $F$ along $\beta$ and the top side $[0,1]\times 1$ of $\mathbb{D}$ attached to the negative side of $F$ along an arc $\beta'$ that is freely isotopic to $h(\beta)$ and intersects $\beta$ minimally. So the remaining two sides of $\mathbb{D}$ are on $\partial X$. We fix an orientation on $\mathbb{D}$ and tilt the disk $\mathbb{D}$ slightly near $F$ so that the positive normal directions of $\mathbb{D}$ and $F$ agree at the branched locus (Figure \ref{fig: branched surface B=FD}). If we view $X = F\times [0,1]/(x, 1)\sim (h(x), 0)$, then the disk $\mathbb{D}$ in $X$ is isotopic to $\beta\times [0,1]/\sim_h$. A key property of the canonical meridian $\mu_0$ defined in \S \ref{subsec: fibered knot} is that: up to isotopy, it is disjoint from the disk $\mathbb{D}\cap \partial X$, the green arcs in Figure \ref{fig: branched surface B=FD}. This follows immediately from the construction of the canonical meridian (see \cite[\S 4]{RobertsSurfacebundle2}).

\begin{figure}[ht]
\centering
\begin{tikzpicture}[scale=0.7]
%\draw [help lines] (0,0) grid (20,7);
%%% Torus on the left
%orange arc
\draw [thick, red] (1.67, 3.1) ..controls (2, 3.1) and (2.6,3.1)  .. (2.8, 3.5);
\draw [thick, red, dashed] (1.3, 3.83) ..controls (2.1, 3.85) and (2.6,3.9)  .. (2.8, 3.5);
\draw [thick, red, -{>[scale=2.5, length=0.7, width=1]}] (2.5, 3.1) -- (2.4, 3.35);
\node [red, above] at (2,3) {\tiny{$\beta$}};
%blue arc
\draw [thick, blue] (1.55, 2.81) .. controls (3,2.5) and (6.3, 3.6) .. (3.55, 3.69);
\draw [thick, blue, dashed] (1.38, 4.1) .. controls (2, 4.1) and (3, 4.2) .. (3.5, 3.7);
\node [blue, below] at (4,3.1) {\tiny{$\beta'$}};
\draw [thick, blue, -{>[scale=2.5, length=0.7, width=1]}] (2.48, 2.88) -- (2.6, 2.6);
%torus body
\draw [thick] (1.5,2.7) to [bend left] (1.5, 4.3);
\draw [thick] (1.5,2.7) to [bend right] (1.5, 4.3);
\draw [thick] (1.5, 2.7) .. controls (2.5,2) and (5.5, 1.5) .. (5.5, 3.5);
\draw [thick] (1.5, 4.3) .. controls (2.5,5) and (5.5, 5.5) .. (5.5, 3.5);
\draw [thick] (2.8, 3.5) to [bend left] (4.1, 3.5);
\draw [thick] (2.8, 3.5) to [bend right] (4.1, 3.5);
%%% product disk
\draw [thick, red] (7.9, 2.3) -- (10.1, 2.3);
\draw [thick, blue] (7.9, 4.7) -- (10.1, 4.7);
\draw [thick, YellowGreen] (7.9,2.3) -- (7.9, 4.7);
\draw [thick,YellowGreen] (10.1, 2.3) -- (10.1, 4.7);
\node at (9, 3.5) {$\mathbb{D}$};
%%% Train track on the boundary
\draw [thick, gray] (12.5, 2.5) to [bend left] (12.5, 4.5);
\draw [thick, gray, dashed] (12.5, 2.5) to [bend right] (12.5, 4.5);
\draw [thick, gray] (18, 2.5) to [bend left] (18, 4.5);
\draw [thick, gray] (18, 2.5) to [bend right] (18, 4.5);
\draw [thick, gray] (12.5, 2.5) -- (18,2.5);
\draw [thick, gray] (12.5, 4.5) -- (18,4.5);
% arcs
\draw [thick,YellowGreen] (17.2, 3.5) .. controls (16.7,3.5) and (17,4.5) .. (16.5, 4.5);
\draw [thick,YellowGreen, dashed] (16.5, 4.5) .. controls (16,4.5) and (16.5,2.5) .. (16,2.5);
\draw [thick,YellowGreen] (16,2.5) .. controls (15.5,2.5) and (16,3.5) .. (15.3, 3.5);
\draw [thick,YellowGreen] (14.5,3.5) .. controls (15,3.5) and (15,4.5) .. (14.5, 4.5);
\draw [thick, YellowGreen,dashed] (14.5, 4.5) .. controls (13.8,4.5) and (14.5,2.5) .. (13.5,2.5);
\draw [thick,YellowGreen] (13.5,2.5) .. controls (13,2.5) and (13.2,3.5) .. (13.7,3.5);
\draw [thick] (12.2, 3.5) -- (17.7, 3.5);
%arrow
\draw [thick, {<[scale=2.5, length=1, width=2]}-] (13, 3.5) -- (13.2, 3.5);
\draw [thick, gray, -{>[scale=2.5, length=1, width=2]}] (17.703,3.7) to  (17.706,3.73);
\node [gray] at (18, 3.5) {\tiny{$\mu_0$}};
\end{tikzpicture}
\caption{An example of the branched surface $B=\langle F;\mathbb{D}\rangle$. The branched locus is $\beta\cup \beta'$, where $\beta'$ is freely isotopic to $h(\beta)$. Since in this example, $\beta$ and $\beta'$ are disjoint, locally the branched surface is homeomorphic to one of the first two models in Figure \ref{fig:branched surface}.  The rightmost figure shows $B\cap \partial X$.}
\label{fig: branched surface B=FD}
\end{figure}
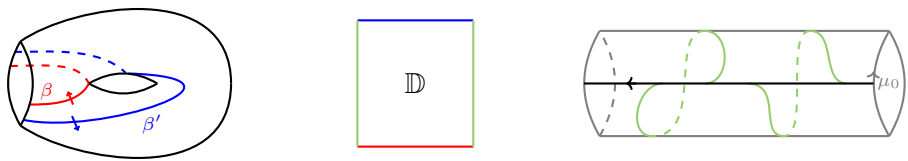

Similarly, given a collection of pairwise disjoint properly embedded arcs $\{\beta_i : i = 1,\cdots, m\}$ in $F$, one can attach $m$ disks $\mathbb{D}_i = [0,1]\times [0,1]$ along $\beta_i$ and $h(\beta_i)$ (up to isotopy) to obtain a branched surface $B=\langle F; \mathbb{D}_1,\cdots ,\mathbb{D}_m \rangle$. 

\end{example}
Branched surfaces provide a useful tool for constructing foliations. We say a foliation $\mathcal{F}$ is {\it (fully) carried} by a branched surface $B$ if after possibly splitting a finite number of leaves, leaves of the foliation can be isotoped to lie in $N(B)$ and intersect every $I$-fiber of $N(B)$ transversely. Intuitively, leaves of the foliation are locally ``parallel'' to $B$ if it is carried by $B$.

\subsection{Computing the Euler class of taut foliations on Dehn fillings of fibered manifolds}  Suppose that $\widehat{\mathcal{F}}$ is a taut foliation on $X(\alpha)$ given in Theorem \ref{thm: foliation fibered knot} and Theorem \ref{thm: foliation 3 braids}, where $X$ is the exterior of a fibered knot in a $\mathbb{Q}$-homology sphere, $\alpha = p\mu_0 + q\lambda$ with $p>0$ and $\mu_0$ is the canonical meridian of $X$, which equals the knot meridian if $X$ is the exterior of a fibered knot in $S^3$.

Since by assumption $\widehat{\mathcal{F}}$ is transverse to the core of the filling solid torus $N$, we assume that $\mathcal{D} = \widehat{\mathcal{F}}|_N$ is the foliation of $N$ by meridional disks. Let $\mathcal{F} = \widehat{\mathcal{F}}|_{X}$. We note that $\mathcal{F}$ is carried by a branched surface $B = \langle F; \mathbb{D}_1, \cdots, \mathbb{D}_{m}\rangle$ as described in Example \ref{ex: branched surface fd}  (\cite[Theorem 4.1, Corollary 4.3, Corollary 4.4]{RobertsSurfacebundle2}  and \cite[\S 3]{krishna18}). 

\begin{proof} [Proof of Theorem \ref{thm: euler class Rachel's foliation}]
As in the proof of Theorem \ref{thm: e=0 necessary sufficient condition}, we will use Lemma \ref{lem:cohomology class zero} with cohomology classes $c$, $c|_X$, $c'$, $c'|_X$ and $c'|_N$ in the lemma replaced by Euler classes $e(\widehat{\mathcal{F}})$, $e(\mathcal{F})$, $e_\sigma(\widehat{\mathcal{F}})$, $e_\sigma(\mathcal{F})$ and $e_\sigma(\mathcal{D})$ respectively. Here, $\sigma$ is again a nowhere vanishing outward pointing section of $T\widehat{\mathcal{F}}$ along $\partial X$. 

Condition (1) of Lemma \ref{lem:cohomology class zero} corresponds to the condition that  $e(\mathcal{F}) = 0$ in Theorem \ref{thm: euler class Rachel's foliation}. Since $X$ is fibered, we have $k=1$. Hence, to complete the proof, it remains to verify that Condition (2) in Lemma \ref{lem:cohomology class zero} is equivalent to $(2g-1) \cdot |q| \equiv  1 \mbox{ \rm{(mod $p$)}}$.

Let $B = \langle F; \mathbb{D}_1, \cdots, \mathbb{D}_{m}\rangle$ be a branched surface that carries $\mathcal{F}$.  Up to isotopy we assume that $\mathcal{F}|_{N(B)}$ is transverse to the $I$-fibers of the normal neighborhood $N(B)$. Since $B$ contains $F$ as a subsurface, we can orient $\mathcal{F}$ so that $T\mathcal{F}|_F$ is homotopic to the tangent field $TF$ of $F$. To see this, we first fix a Riemannian metric $g$ on $X$. Let $v_0$ and $v_1$ be the positive normal vector fields of $T\mathcal{F}|_F$ and $TF$ respectively. We assume that $v_1$ is tangent to the $I$-fibers of $N(F) = F\times I$, and both $v_0$ and $v_1$ are tangent to $\partial X$.   The condition that $T\mathcal{F}|_{N(B)}$ is transverse to the $I$-fibers of $N(B)$ is equivalent to that $v_0\cdot v_1\neq 0$, where the dot product $\cdot$ is taken pointwise. We assume that $v_0\cdot v_1 > 0$ over $F$. Otherwise, we reverse the orientation of $\widehat{\mathcal{F}}$ (and hence that of $\mathcal{F}$) and replace $v_0$ by $-v_0$. Let $v_t=tv_0+(1-t)v_1$ for $t\in [0,1]$. Note that $v_t$ is tangent to $\partial X$ for all $t$. Moreover, $v_t\cdot v_1=t v_0\cdot v_1+(1-t)\Vert v_1 \Vert$ are positive and hence nonzero.  So the desired homotopy $\xi_t$ is given by the plane fields that are normal to $v_t$ for $t\in [0,1]$.  

Let $\sigma$ be an outward pointing nowhere vanishing vector field that is tangent to $\xi_0= T\mathcal{F}|_F$ and $\xi_1 =  TF$ along $\partial F$. So $\sigma$ is tangent to all $\xi_t$ for $t\in [0,1]$. Therefore, we obtain a relative Euler class of $\xi_t$ as a plane field over $F\times [0,1]$ in $H^2(F\times [0,1], \partial F\times [0,1])$ associated with the section $\sigma\times {\rm id}$ over $\partial F\times [0,1]$, whose projections to $F\times 0$ and $F\times 1$ are $e_\sigma(T\mathcal{F}|_F)$ and $e_\sigma(TF)$ respectively. 
Hence, $e_\sigma(T\mathcal{F}|_F)= e_\sigma(TF)$. By the Poincare-Hopf Theorem, $e_\sigma(TF)([F]) = \chi(F)$ and we have 
$$a =e_\sigma(T\mathcal{F})([F])=e_\sigma(T\mathcal{F}|_F)([F]) =  \chi(F) = 1-2g.$$ 

It remains to determine the value of $b = e_\sigma(\mathcal{D})([D])$. To do so, we need to compare the orientation of the foliation $\mathcal{D}$ and the given orientation on the meridian disk $D$ of $N$. This is done by comparing their induced orientations on the simple closed curves  on $\partial X$.

Recall that the meridian disk $D$ is oriented so that $f(\partial D) = p\mu_0 + q\lambda$ with $p>0$, where $f: \partial N\rightarrow \partial X$ is the gluing map (See \S \ref{subsec: notation and orientation}).  Let $\alpha'$ be a component of $\mathcal{D}\cap \partial X$ contained in $N(B)\cap \partial X$ with the induced boundary orientation from $\mathcal{D}$. So  $\alpha' = \pm (p\mu_0 + q \lambda)$ and $-\alpha'$ is a component of $\mathcal{F}\cap \partial X$ with the induced boundary orientation from $\mathcal{F}$. 

Since $\alpha'$ is contained in $N(B)$ and the canonical meridian $\mu_0$ can be isotoped to be disjoint from $\mathbb{D}_i\cap \partial X$ for $i=1,\cdots, m$ (Example \ref{ex: branched surface fd}), we have $\mu_0$ intersects only  the portion of $\alpha'$ in $N(\partial F)$, the normal neighborhood of the black line in the rightmost picture of Figure \ref{fig: branched surface B=FD}). Because we oriented $\mathcal{F}$ so that leaves of $\mathcal{F}$ intersects the $I$-fibers of $N(F)$ positive, it follows that $-\alpha'$ intersects $\mu_0$ positively. 

Hence, $- \alpha' = p\mu_0 + q \lambda$ if $p/q>0$ and $-\alpha' = - p\mu_0 - q \lambda$ if $p/q<0$. Therefore,

(1) If $p/q > 0$, then $\alpha' = -(p\mu_0 + q \lambda)$. So $\mathcal{D}$ has the opposite orientation with $D$ and hence $e_\sigma(\mathcal{D})([D]) = -1$. By Lemma \ref{lem:cohomology class zero}, we have $e(\widehat{\mathcal{F}}) = 0$ if and only if $e(\mathcal{F}) = 0$ and  
\begin{displaymath}
 \chi(F) \cdot q \equiv - 1 \mbox{ \rm{(mod $p$)}}. 
\end{displaymath}
Because $p/q > 0$ and $p>0$ implies that $q>0$, the above identity is equivalent to $|\chi(F) \cdot q| =  1 \mbox{ \rm{(mod $p$)}}$

(2) If $p/q < 0$, then $\alpha' = p\mu_0 + q \lambda$. So $\mathcal{D}$ has the same orientation with $D$ and hence $e_\sigma(\mathcal{D})([D]) = 1$. By Lemma \ref{lem:cohomology class zero}, we have $e(\widehat{\mathcal{F}}) = 0$ if and only if $e(\mathcal{F}) = 0$ and  
\begin{displaymath}
 \chi(F) \cdot q \equiv  1 \mbox{ \rm{(mod $p$)}}  
 \end{displaymath}
Because $p/q < 0$ and $p>0$ implies that $q<0$, the above identity is again equivalent to $|\chi(F) \cdot q| =  1 \mbox{ \rm{(mod $p$)}}$.
\end{proof}

\subsection{Computing the slopes in $\mathcal{L}_g$}
\label{sec: computation of the slopes}

For each $g>0$, we denote the set of slopes satisfying Equation $(\ref{equ: fiber knot slope equation})$  in Theorem \ref{thm: euler class Rachel's foliation}  by 
$$\mathcal{L}_g = \{p/q \in \mathbb{Q}\setminus\{0\} : (2g-1) |q|  \equiv 1 \mbox{ \rm{(mod $p$)}}\}.$$
 
\subsubsection{Genus one knots}
\label{subsec: genus one}
When $g=1$, the above equation reduces to 
\begin{equation}
 |q|\equiv 1 \mbox{ (mod $p$)},
 \label{equ: branched surface quality g=1}
\end{equation}
 which is equivalent to $|q| = 1+np$ for some $n\geq 0$. Hence $q = \pm (1+np)$ with $n\geq 0$ and $p\geq 1$.

\begin{lemma}
 Let \begin{equation}
 \mathcal{L}_1 = \left\{\pm \frac{p}{1+np} : p\geq 1 \text{ and } n \geq 0\right\}.
 \label{equ: slope branched surface g=1}
 \end{equation}
Then $\mathcal{L}_1$ is the set of slopes for which Equality (\ref{equ: fiber knot slope equation}) holds and $\mathbb{Z}\setminus 0 \subset \mathcal{L}_1$.  
\label{lem: slope branched surface g=1}
\end{lemma}
\begin{proof}
 One obtains all nonzero integer slopes by taking $n=0$.
\end{proof}

\subsubsection{Higher genus knots}  
\label{subsec: higher genus}
\begin{lemma}
\label{lem: slope branched surface g>1}
The set $\mathcal{L}_g$  is bounded between $\pm (2g-1)$. Let $\alpha_k = (2g-1) - \frac{1}{k+1}$, $k\in \mathbb{N}$. Then $\pm\alpha_k$  are in $\mathcal{L}_g$ for all $k\in \mathbb{N}$. 
\end{lemma}

\begin{proof}
Given a slope $p/q\in \mathcal{L}_g$, $p>0$,   it satisfies 
 \begin{displaymath}
  (2g-1)|q|=1 + np \text{ for some } n\in \mathbb{N}.
 \end{displaymath}
Because $g>1$, $n$ cannot be $0$ in the above equation. Hence $p/q\in \mathcal{L}_g$  is equivalent to 
 \begin{displaymath}
   \frac{p}{|q|} = \frac{1}{n} \left((2g-1) -\frac{1}{|q|}\right). 
 \end{displaymath} 
 
 It follows that the set $\mathcal{L}_g$ is bounded between $\pm (2g-1)$. One obtains $\pm \alpha_k$ by taking  $n=1$ and $q=\pm (k+1)$ for $k\in \mathbb{N}$. 
\end{proof}

\subsection{The proof of results on group left-orderability}
\label{subsec: lo fibered manifolds}
Suppose that $X$ is the exterior of  a knot in a $\mathbb{Z}$-homology sphere, then $H^2(X) = 0$ and the condition that $e(\mathcal{F}) =0$ in Theorem \ref{thm: euler class Rachel's foliation} always holds. Since Corollary \ref{cor: 3 braids left-orderability} only concerns knots in $S^3$, it follows from Theorem \ref{thm: euler class Rachel's foliation} and Lemma \ref{lem: slope branched surface g>1}. 

Next we prove Corollary \ref{cor: left-orderable group g=1 fibered} and an analogous result for fibered knots of higher genus.

\begin{proof}[Proof of Corollary \ref{cor: left-orderable group g=1 fibered}]
When $\alpha = \lambda$ is the zero slope, the closed $3$-manifold is a surface bundle over $S^1$. So it is irreducible and $b_1(X(\alpha)) > 0$, which shows that $\pi_1(X(\alpha))$ is LO by \cite[Corollary 3.4]{BRW}. 

Now suppose that $\alpha \neq \lambda$. Let $\widehat{\mathcal{F}}$ denote a taut foliation on $X(\alpha)$ given in Theorem \ref{thm: foliation fibered knot}. We will show that the Euler class of $e(\widehat{\mathcal{F}})$ is zero and hence by Theorem \ref{thm: universal circle}, we have $\pi_1(X(\alpha))$ is LO. 

Let $\alpha = p\mu_0 + q\lambda \in \mathcal{L}_1$ with $p>0$. Since $\alpha \in \mathcal{L}_1$, we already have the condition that $(2g-1)|q| \equiv 1 \mbox{ \rm{(mod $p$)}}$ satisfied and it remains for us to show that $e(\mathcal{F}) = 0$ in $H^2(X)$, where as before $\mathcal{F} = \widehat{\mathcal{F}}|_X$. This is done in \cite[\S 10]{BH19}. We outline the proof here (See the proof of \cite[Theorem 1.9]{BH19}).  

In the proof of Theorem \ref{thm: euler class Rachel's foliation}, we have shown that the plane field $T\mathcal{F}$ is homotopic to $TF$ over $F$. In fact, since the taut foliation $\mathcal{F}$ on $X$ is carried by a branched surface of the form $B=\langle F; \mathbb{D} \rangle$, the entire tangent plane field $T\mathcal{F}$ is homotopic to the tangent plane field of the fibration on $X$, which is the plane field over $X$ tangent to the fibers of $X$ (\cite[Lemma 4.4]{HKMII}).

Let $h: F\rightarrow F$ be the monodromy of the fibered knot with $h_{\partial F} = {\rm id}$ and $M=(F,h)$ denote the closed manifold given by the open book $(F,h)$.  Let $\tilde{\xi}$ denote the contact structure supported by $(F,h)$. In \cite[Theorem 1.9]{BH19}, it is shown that  the Euler class of the contact structure $e(\tilde{\xi}) =0$ in $H^2(M)$. 

Because $\tilde{\xi}$ is supported by the open book $(F, h)$, by definition its restriction $\tilde{\xi}|_X$ is also homotopic to the tangent plane field of the fibration on $X$. Therefore,  $e(\mathcal{F}) = e(\tilde{\xi}|_X) =0$ in $H^2(X)$. 
\end{proof}

A key fact from the condition that the genus of $F$ is one used implicitly in the above proof is that the mapping class group ${\rm MCG}(F)$ is isomorphic to the group of $3$-strand braid $B_3\cong {\rm MCG}(D_3)$, where $D_3$ is the disk with three marked points. The isomorphism from ${\rm MCG}(D_3)$ to ${\rm MCG}(F)$  is given by lifting elements in ${\rm MCG}(D_3)$ to homeomorphisms of $F$ through the twofold branched cover $F\rightarrow D_3$. The fact that $h$ is a lift from a mapping class of $D_3$ is what allows one to conclude that the Euler class of the contact structure supported by $(F,h)$ is zero in the proof of  \cite[Theorem 1.9]{BH19}.

In general, mapping classes of a surface $F$ of genus $g$ with connected boundary that can be realized as lifts of elements in ${\rm MCG}(D_{2g+1})$ are called {\it symmetric homeomorphisms}. When $g>1$,  they form a proper subgroup of ${\rm MCG}(F)$, called the {\it symmetric mapping class group} (\cite{BirH}; also see \cite[\S 9.4]{FM}).  In \cite[Theorem 1.9]{BH19}, it was shown that  as long as the monodromy $h: F\rightarrow F$ is a symmetric homeomorphism, the Euler class of the contact structure supported by $(F, h)$ is zero. Hence, the argument in the proof of Corollary \ref{cor: left-orderable group g=1 fibered} leads us to the following.

\begin{corollary}
 Let $K$ be a hyperbolic fibered knot in a $\mathbb{Q}$-homology sphere with the degeneracy slope $\gamma$. Suppose that the monodromy $h$ of $K$ is a symmetric homeomorphism. Then with respect to the canonical meridian $\mu_0$, we have 
 \begin{enumerate}
  \item If $\gamma = \mu_0$, then given any $\alpha \in \mathcal{L}_g$, there exists a co-orientable taut foliation on $X(\alpha)$ of zero Euler class and $\pi_1(X(\alpha))$ is LO.   
 \item If $\gamma$ is a positive slope, then given any $\alpha \in \mathcal{L}_g\cap (-\infty, 1)$, there exists a co-orientable taut foliation on $X(\alpha)$ of zero Euler class and hence $\pi_1(X(\alpha))$ is LO.   
 \item If $\gamma$ is a negative slope, then given any $\alpha \in \mathcal{L}_g \cap (-1, \infty)$, there exists a co-orientable taut foliation on $X(\alpha)$ of zero Euler class and hence $\pi_1(X(\alpha))$ is LO.  
 \end{enumerate}
 \label{cor: left-orderable group g>1 fibered}
\end{corollary}

\section{Euler class and restrictions of filling slopes}
\label{sec: euler class slopes}
Given $X$ a $\mathbb{Q}$-homology solid torus,  let $\mathcal{S}_X$ be the set of slope defined as follows:  for any slope $
\alpha\in\mathcal{S}_X$, there exists a taut foliation on $X(\alpha)$ transverse to the core of the Dehn filling solid torus whose Euler class is zero.
\begin{enumerate}
 \item By Theorem \ref{thm: universal circle}, given any $\alpha\in \mathcal{S}_X$, we have $\pi_1(X(\alpha))$ is left-orderable.
 \item Given a taut foliation $\widehat{\mathcal{F}}$ on $X(\alpha)$ for which the core of the filling solid torus is a closed transversal, if $\alpha\notin \mathcal{S}_X$,  then the Euler class $e(\widehat{\mathcal{F}})\neq 0$. 
\end{enumerate}

We fix a meridian $\mu$ and identify the set $\mathcal{S}_X$ with a subset of $\mathbb{Q}\cup\{\frac{1}{0}\}$. We denote it by  $\mathcal{S}_{X, \mu} \subset\mathbb{Q}\cup\{\frac{1}{0}\}$.  The aim of this section is to gain a better understanding of the distribution of the set $\mathcal{S}_{X,\mu}$ in $\mathbb{R}\cup \{\frac{1}{0}\}$. We prove Theorem \ref{thm: restriction on SX null homologous} in this section. 

\subsection{The Thurston norm of the relative Euler class of a taut foliation}
\label{subsec: relative class taut foliation} 
Let $\widehat{\mathcal{F}}$ be a taut foliation on $X(\alpha)$ transverse to the core of the filling solid torus and such that $\mathcal{F} = \widehat{\mathcal{F}}|_X$ intersects $\partial X$ in simple closed curves of slope $\alpha$. Since $\mathcal{F}$ is taut, it is a well-known result that the Thurston dual norm of the relative class $e_\sigma(\mathcal{F}) \in H^2(X,\partial X)$ is at most $1$ \cite[\S 3, Corollary 1]{Thurston86}.  This follows from the fact that any properly embedded connected incompressible and boundary incompressible surface in a tautly foliated $3$-manifold can be isotoped to be either a leaf of the foliation or transverse to the foliation everywhere except for at a finite number of saddle points \cite[\S 3, Theorem 4]{Thurston86} (also see \cite[Theorem 3.7]{Gab00}).  We state Thurston's result in the form that is the most convenient for us in Theorem \ref{thm: norm bound} below. 

\begin{theorem}[Thurston]
Let $X$ be an oriented $\mathbb{Q}$-homology solid torus and $\mathcal{F}$ be a co-oriented taut foliation on $X$ that is transverse to $\partial X$.  Let $F$ be an oriented connected incompressible surface properly embedded in $X$ such that $[F]$ represents a generator of $H_2(X,\partial X)$. Assume that each boundary component of $F$ is either tangent or transverse to $\mathcal{F}|_{\partial X}$. Then the relative Euler class $e_\sigma(\mathcal{F}) \in H^2(X,\partial X)$, associated to a nowhere vanishing outward-pointing (or inward-pointing) section $\sigma: \partial X\rightarrow T\mathcal{F}|_{\partial X}$ satisfies:
 \begin{itemize}
 \item [(i)] $|e_\sigma(\mathcal{F})([F])|\leq |\chi(F)|$. 
 \item [(ii)] $e_\sigma(\mathcal{F})([F]) \equiv \chi(F) \mbox{ \rm{(mod $2$})}$.
\end{itemize}
\label{thm: norm bound}
\end{theorem}

\begin{proof}
We may assume that $\sigma$ is tangent to $F$ along $\partial F$. If $F$ is isotopic to a leaf of the foliation $\mathcal{F}$, then up to isotopy, $T\mathcal{F}|_F= \pm TF$. By the Poincare-Hopf Theorem, we have  $e_\sigma([F])=\pm \chi(F)$,  where the sign is $+$ if the orientations on $F$ and $\mathcal{F}$ are the same and $-$ otherwise.

Now suppose that $F$ cannot be isotoped to a leaf.  Since $\mathcal{F}$ is taut, $X$ is irreducible \cite{Nov65,Ros68}. So the incompressibility of $F$ implies that $F$ is also boundary incompressible \cite[Lemma 1.10]{Hatcher} and hence $F$ can be isotoped to be transverse to $\mathcal{F}$ except at a finite number of saddle (index $-1$) tangential points in the interior of $F$ \cite{Thurston86, Gab00}. Let $\mathcal{L}=\mathcal{F}\cap F$ be the induced singular foliation on $F$. At each singular point $p$ of $\mathcal{L}$, we have $T_p\mathcal{F}=\pm T_pF$. 
Let $e_+$ (resp. $e_-$) be the total number of singular points on $\mathcal{L}$, where the orientation of $TF$ is the same with (resp. opposite to) the orientation on $T\mathcal{F}$. Hence  $\chi(F)=-e_+-e_-$ and $ e_\sigma([F])=-e_++e_-$ by the Poincar\'e-Hopf Thoerem and the definition of Euler class. It follows that $|e_{\sigma}([F])|\leq |\chi(F)|$ and that $e_\sigma([F])$ has the same parity with $\chi(F)$. 
\end{proof}

\subsection{Restrictions on the Dehn filling slopes}
Let $\mathsf{x}$ be the Thurston norm (\cite{Thurston86}) of a generator $c$ of $H_2(X,\partial X)\cong \mathbb{Z}$, which is defined to be 
\begin{displaymath}
\min \left\{\max\{0, -\chi(F)\} \mid F \text{ is oriented properly embedded in } X \text{ and } [F] = c\right \}.
\end{displaymath}

Theorem \ref{thm: restriction on SX null homologous} is the special case of Theorem  \ref{thm: restriction on SX} below  with 
$k=1$ and $\mathsf{x} = 2g-1$, $g>0$. 

\begin{theorem}
Let $X$ be a $\mathbb{Q}$-homology solid torus, $k\geq 1$ the order of the longitude of $X$ and $\mathsf{x}$ the Thurston norm of a generator of $H_2(X,\partial X)$. Assume that $\mathsf{x}\neq 0$. Then given any meridian $\mu$, we have the following: 
\begin{enumerate}
 \item Outside $\left(-\frac{\mathsf{x}}{k} - 1, \frac{\mathsf{x}}{k} + 1\right)$, the set $\mathcal{S}_{X, \mu}$ only contains $\mu$ and the integer slopes. That is, 
 $$\mathcal{S}_{X, \mu}\setminus \left(-\frac{\mathsf{x}}{k} - 1, \frac{\mathsf{x}}{k} + 1\right) \subseteq \mathbb{Z}\cup \left\{\frac{1}{0}\right\}.$$
 \item The set $\mathcal{S}_{X,\mu}$ is nowhere dense in $\mathbb{R}\cup\{\frac{1}{0}\}\cong \mathbb{R}P^1$. Particularly, it is nowhere dense within $\left(-\frac{\mathsf{x}}{k} - 1, \frac{\mathsf{x}}{k} + 1\right)$.
 \end{enumerate}
\label{thm: restriction on SX}
\end{theorem}

There is a useful corollary that follows immediately from Theorem \ref{thm: restriction on SX}. It gives a simple criterion for a slope to be in $\mathcal{S}_X$.

\begin{corollary}
Let $X$ be a $\mathbb{Q}$-homology solid torus, $k\geq 1$ be the order of the rational longitude and $\mathsf{x}$ the Thurston norm of a generator of $H_2(X, \partial X)$. Fix a meridian $\mu$. Given   
 a co-oriented taut foliation $\widehat{\mathcal{F}}$  on $X(p/q)$ that is transverse to the core of the filling solid torus,  if  $|\frac{p}{p+nq}|>\frac{\mathsf{x}}{k}+1$ for some $n\in \mathbb{Z}$ and $\frac{p}{p+nq} \notin\mathbb{Z}\cup \{\frac{1}{0}\}$, then the Euler class $e(\widehat{\mathcal{F}})\neq 0$.
\label{cor: obstruction k>1}
\end{corollary}

\begin{proof}
Given any $n\in \mathbb{Z}$, let $\mu' = \mu -n \lambda$ be another meridian. Then $$\alpha = p\mu + q\lambda = p\mu' + (np + q)\lambda.$$ Hence $\alpha = \frac{p}{q + np}$ with respect to $\mu'$. By Theorem \ref{thm: restriction on SX}(1), since $|\frac{p}{p+nq}|>\frac{\mathsf{x}}{k}+1$ and $\frac{p}{p+nq} \notin\mathbb{Z}\cup \{\frac{1}{0}\}$, we have $\frac{p}{q + np}\notin \mathcal{S}_{X,\mu'}$. Therefore, $e(\widehat{\mathcal{F}})\neq 0$.  
\end{proof}

For instance, let $X$ be the exterior of a genus one knot in an $\mathbb{Z}$-homology sphere. So we have $k=1$ and $\mathsf{x}=1$. Consider slope $\frac{5}{12}$. Since $$\frac{5}{12 - 2\times 5} = \frac{5}{2} >  |\chi(F)| + 1 =2,$$ by Corollary \ref{cor: obstruction k>1} we can conclude that $\frac{5}{12}$ is not in $\mathcal{S}_X$. Consequently, the Euler class of any taut foliation on $X(\frac{5}{12})$ that is transverse to the core of the filling solid torus is nonzero. 

\begin{remark}\label{rem: theorem 4.2 weaker form}
By the definition of the Thurston norm, both Theorem \ref{thm: restriction on SX} and Corollary \ref{cor: obstruction k>1} hold if we replace $\mathsf{x}$ by $-\chi(F)$, where $F$ is any connected, oriented, properly embedded incompressible surface in $X$ that represents a generator of $H_2(X,\partial X)$. Namely, $F$ doesn't need to be norm-minimizing, though one achieves the best bound if it is. 
\end{remark}

\begin{proof}[Proof of Theorem \ref{thm: restriction on SX}(1)]
Let $\alpha = p\mu+q\lambda$ in $\mathcal{S}_{X}$ with $p\geq 0$.  Since $p/q\in \mathcal{S}_{X,\mu}$ there exists a co-orientable taut foliation, denoted by $\widehat{\mathcal{F}}$ on $X(p/q)=X\cup_f N$ that is transverse to the core of $N$ and $e(\widehat{\mathcal{F}}) = 0$. Shrinking $N$ if necessary, we may assume that $\mathcal{F} =\widehat{\mathcal{F}}|_X$ intersects $\partial X$ transversely in simple closed curves of slope $p/q$  and $\mathcal{D}= \widehat{\mathcal{F}}|_N$ is the foliation of $N$ by meridional disks.  We show that if $p/q\notin \mathbb{Z}\cup \{\frac{1}{0}\}$, then $|p/q|\leq \frac{\mathsf{x}}{k}+1$. This proves Part (1) of the theorem.

Let $F$ denote a connected, oriented, properly embedded norm-minimizing surface which represents a generator in $H_2(X,\partial X)$. Since $p/q\notin \mathbb{Z}$, we have $p\neq 0$. We orient $\widehat{\mathcal{F}}$ as in Theorem \ref{thm: e=0 necessary sufficient condition}. Then by Theorem \ref{thm: e=0 necessary sufficient condition}(2) we have $e(\widehat{\mathcal{F}}) = 0$ implies that $(aq)/k  \equiv 1 \mbox{ \rm{(mod $p$)}}$, where $a = e_\sigma(\mathcal{F})([F])$ is an integer divisible by $k$. We let $a'=a/k$.  So the equation $(aq)/k  \equiv 1 \mbox{ \rm{(mod $p$)}}$ is equivalent to $a'q= 1 + sp$ for some $s\in \mathbb{Z}$.
Also because $p/q\notin \mathbb{Z}\cup \{\frac{1}{0}\}$, we have $s\neq 0$ and $q\neq 0$. Hence, $a'q=1 +sp$ is also equivalent to 
\begin{equation}
 \frac{p}{q} = \frac{1}{s} (a' - \frac{1}{q}).
 \label{equ: k>1 M_g}
\end{equation}

By Theorem \ref{thm: norm bound},  $a'=\frac{a}{k}$ is between $\pm \chi(F)/k$. Therefore, 
\begin{displaymath}
 \left|\frac{p}{q}\right| \leq \left|\frac{1}{s} (a' - \frac{1}{q})\right| \leq |a'|+1 \leq |\chi(F)|/k+1 = \frac{\mathsf{x}}{k} +1. 
\end{displaymath}
This completes the proof of Theorem \ref{thm: restriction on SX} Part (1).
\end{proof}

To prove Theorem \ref{thm: restriction on SX}(2), we construct a slightly larger set denoted by $\mathcal{M}_g$ if $k=1$, and $\mathcal{M}_\mathsf{x}$ in general. This set is highly symmetric and easy to visualize, which helps us gain a better intuition of the set $\mathcal{S}_{X}$. We will show that the set $\mathcal{M}_{\mathsf{x}}$ is nowhere dense in $\mathbb{R}\cup\{\frac{1}{0}\}$. Since $\mathcal{S}_{X}$ is a subset of it, it is also nowhere dense.

We first consider the case when $k=1$. This is equivalent to $X$ being the exterior of a null-homologous knot $K$ in a $\mathbb{Q}$-homology sphere. 

\begin{proposition}
Let $X$ be the exterior of a null-homologous knot of genus  $g>0$ in an oriented $\mathbb{Q}$-homology sphere. Let $\mathcal{M}_g$ be the subset of $\mathbb{Q}$ that is uniquely determined by the following properties:  
\begin{enumerate}
\item $ \mathcal{M}_g = \mathbb{Z} \cup (\cup_{n\geq 1} \, \mathcal{M}_{g ,n})$.
 \item $\mathcal{M}_{g,1}=\{m+\frac{1}{s} \mid s\in \mathbb{Z}\setminus \{0\}, \text{ and $m$ is an odd number between $\pm(2g-1)$}\}$.
 \item $\mathcal{M}_{g,n}=\frac{1}{n} \mathcal{M}_{g,1}$ for $n\geq 1$.
\end{enumerate}
Then $\mathcal{S}_{X,\mu}$ is a proper subset of $\mathcal{M}_g\cup \{\frac{1}{0}\}$ for any given meridian $\mu$. 
\label{prop: slope obstruction M_g} 
\end{proposition}

An example of $\mathcal{M}_g$ when $g=2$ is depicted in the Figure \ref{fig: M_g} below.

\begin{figure}[ht]
\centering 
\begin{tikzpicture}[scale = 0.55]
%\draw [help lines] (-12,-4) grid (10,8);
% M_{g,1} 
% label
\node [left] at (-9.5, 5) {\small $\mathcal{M}_{g,1}:$};
\foreach \value in {0,1,2,3,4,-1,-2,-3,-4}
	\node [below] at (2*\value,5) {\small $\value$};
% M_{g,1}
\foreach \a in {-1, 1, -3, 3}
	\foreach \s in {1,2,3,...,30}
	{
	\draw[fill=black] (2*\a+2/\s, 5) circle (0.04);
	\draw[fill=black] (2*\a-2/\s, 5) circle (0.04);
	}
%% This fill white is because the points are too dense near \pm 1 and \pm 3
\foreach \a in {-1, 1, -3, 3}
	\draw[fill = white]  (2*\a, 5) circle (0.04);
	
%% M_{g,2} label
\node [left] at (-9.5, 3) {\small $\mathcal{M}_{g,2}:$};
\foreach \value in {2,1,0,-1,-2}
	\node [below] at (2*\value, 3) {\small $\value$};
% M_{g,2} 
\foreach \a in {-1, 1, -3, 3}
	\foreach \s in {1,2,3,...,15}
		{\draw[fill=black] (\a+1/\s, 3) circle (0.04);
		 \draw[fill=black] (\a-1/\s, 3) circle (0.04);
		}
\foreach \a in {-1, 3, -3, 1} 
	\draw [fill = white] (\a, 3) circle (0.04);
% M_{g,4} label
\node [left] at (-9.5, 1) {\small $\mathcal{M}_{g,4}:$};
\foreach \value in {1,0,-1}
	\node [below] at (2*\value, 1) {\small $\value$};
%% M_{g,4}
\foreach \a in {-1, 1, -3, 3}
	\foreach \s in {1,2,3,...,7}
		{\draw[fill=black] (\a/2+0.5/\s, 1) circle (0.04);
		 \draw[fill=black] (\a/2-0.5/\s, 1) circle (0.04);
		}
\foreach \a in {-1, 1, -3, 3}
	{
	\draw[fill =white]  (\a/2, 1) circle (0.04);
	\draw (\a/2, 1) circle (0.04);
	}
% M_g the whole set
\node [left] at (-9.5, -1) {\small $\mathcal{M}_g:$};
\draw [->] (-9 , -1) -- (8.5, -1);
\foreach \Point/\PointLabel in {(-8, -1)/-4, (-6, -1)/-3, (-4, -1)/-2,(-2, -1)/-1,(0, -1)/0,(2, -1)/1, (4,-1)/2, (6,-1)/3, (8,-1)/4}	\draw[fill=black] \Point circle (0.04) node[below] {\small $\PointLabel$};
\foreach \a in {-3, -1, 1, 3}
	\foreach \s in {1, 1/2,1/3,1/4,1/5,1/6,1/7,1/8,1/9, 1/10, 1/11, 1/12, 1/13, 1/14,1/15, 1/16, 1/17, 1/18, 1/19, 1/20, 1/21, 1/22, 1/23}
		\foreach \n in {1,2}
		{
		 \draw[fill=black] (2*\a/\n+2*\s/\n, -1) circle (0.03);
		 \draw[fill=black] (2*\a/\n-2*\s/\n, -1) circle (0.03);
		}
\foreach \a in {-3, -1, 1, 3}
	\foreach \s in {1, 1/2,1/3,1/4,1/5,1/6,1/7,1/8,1/9, 1/10}
		\foreach \n in {2,3,4}
		{
		 \draw[fill=black] (2*\a/\n+2*\s/\n, -1) circle (0.03);
		 \draw[fill=black] (2*\a/\n-2*\s/\n, -1) circle (0.03);
		}
\foreach \a in {-3, -1, 1, 3}
	\foreach \s in {1, 1/2,1/3,1/4,1/5}
		\foreach \n in {5,6,...,10,11,12}
		{
		 \draw[fill=black] (2*\a/\n+2*\s/\n, -1) circle (0.03);
		 \draw[fill=black] (2*\a/\n-2*\s/\n, -1) circle (0.03);
		}
\end{tikzpicture}
\caption{\small Points indicate slopes in $\mathcal{M}_g = \mathbb{Z}\cup (\cup_{n\geq 1} \, \mathcal{M}_{g ,n})$ as defined in Proposition \ref{prop: slope obstruction M_g} with $g=2$. We plotted values until points around limit points became indistinguishable. One can see that $\mathcal{M}_{g,1}$ is contained in the interval $[-2g,2g] = [-4, 4]$ and consists of $8$ sequences: $\{m \pm \frac{1}{k} : m = \pm 1, \pm 3 \text{ and } k\in \mathbb{N}\}$ as described in Proposition \ref{prop: slope obstruction M_g}(2). By Proposition \ref{prop: slope obstruction M_g}(3), each $\mathcal{M}_{g,n}$ for $n> 1$ can be obtained by shrinking $\mathcal{M}_{g,1}$ proportionally. Though the points become denser and denser as they get closer to $0$,  the set $\mathcal{M}_g$ is in fact nowhere dense by Proposition \ref{prop: M_g is nowhere dense}.}
\label{fig: M_g}
\end{figure}
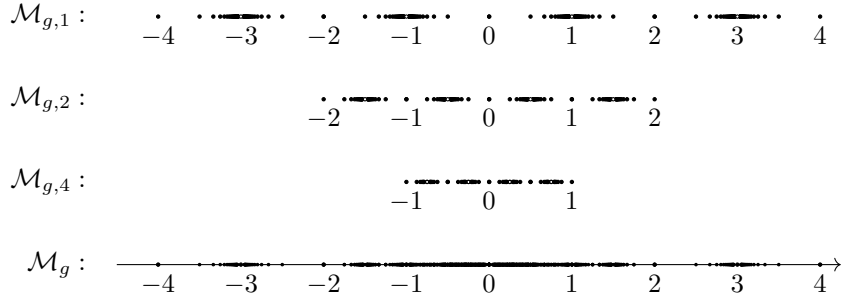

\begin{proof}[Proof of Proposition \ref{prop: slope obstruction M_g}]
Let $\mu$ be a meridian. Suppose that $p/q$ is a slope in $\mathcal{S}_{X,\mu}$. Then by the definition of $\mathcal{S}_{X,\mu}$, there exists a co-orientable taut foliation $\widehat{\mathcal{F}}$ with Euler class zero on $X(p/q)$ that is transverse to the core of the filling solid torus.  Shrinking $N$ if necessary, we assume that $\mathcal{D}= \widehat{\mathcal{F}}|_N$ is the foliation of $N$ by meridional disks. We also orient $\widehat{\mathcal{F}}$ as in Theorem \ref{thm: e=0 necessary sufficient condition}. We want to show that  if $p/q \neq \frac{1}{0}$ (so $q\neq 0$), then $p/q$ is in $\mathcal{M}_g$. 

 Let $F$ be an oriented norm-minimizing Seifert surface.  Since the Euler class $e(\widehat{\mathcal{F}})=0$, by Theorem \ref{thm: e=0 necessary sufficient condition}, we have $aq=1$ (mod $p$), where $a=e_\sigma(\mathcal{F})([F])$.  So there exists an integer $n_0\in \mathbb{Z}$, such that $aq = 1 + n_0p$. If $n_0=0$, then $aq = 1$, which implies that $q = \pm 1$. Hence $p/q \in \mathbb{Z} \subset \mathcal{M}_g$. Suppose that $n_0\neq 0$. Since $q\neq 0$, we have $aq = 1 + n_0p$ is equivalent to
\begin{align*}
 \frac{p}{q} & = \frac{1}{n_0} (a - \frac{1}{q}) = \pm \frac{1}{|n_0|} ( a - \frac{1}{q})  \\
  & = \frac{1}{|n_0|} (\pm  a -\frac{ \pm  1}{q}).
\end{align*}
By Theorem \ref{thm: norm bound}, $a=e_\sigma(\mathcal{F})([F])$ is an odd integer bounded between $\pm (2g-1)$. Therefore, the slope $\frac{p}{q}$ lies in $\frac{1}{n}\mathcal{M}_{g,1} = \mathcal{M}_{g, n} \subset \mathcal{M}_g$, where $n=|n_0|$.
\end{proof}

\begin{proposition}
$\mathcal{M}_g$ is nowhere dense in $\mathbb{R}$. \label{prop: M_g is nowhere dense}
\end{proposition}

\begin{proof}
We show that the closure of $\mathcal{M}_g$ in $\mathbb{R}$, denoted by $\overline{\mathcal{M}}_g$ has empty interior. 

Let $(a, b)$ be an open interval in $\mathbb{R}$. We claim that $(a, b)\nsubseteq \overline{\mathcal{M}}_g$ (i.e. $(a, b)\cap \overline{M}_g\neq (a,b)$).  By taking a sub-interval of $(a, b)$ if necessary, we may assume that $0\notin [a, b]$. We consider the case that $(a,b)\subset (-\infty, 0)$. When $(a,b) \subset (0, \infty)$, the argument is similar. 

Note that since $\mathcal{M}_{g,1}$ is bounded between $\pm 2g$, the set $\mathcal{M}_{g,n}$ is bounded between $\pm \frac{2g}{n}$ for $n\geq 1$.  Let $n_0>0$ be an integer sufficiently large so that the closure of $\cup_{n\geq n_0} \mathcal{M}_{g,n}$ is disjoint from $(a, b)$. Hence 
\begin{displaymath}
   (a,b)\cap \overline{\mathcal{M}_g} = (a,b) \cap  \overline{\cup_{n=1}^{n_0}\mathcal{M}_{g,n}}.
\end{displaymath}
It is easy to see that $\mathcal{M}_{g,1}$ is nowhere dense in $\mathbb{R}$ and hence each $\mathcal{M}_{g,n} = \frac{1}{n}\mathcal{M}_{g,1}$ is nowhere dense.  A finite union of nowhere dense sets is again nowhere dense, so we have $\cup_{n=1}^{n_0}\mathcal{M}_{g,n}$ is nowhere dense. Therefore, $(a, b)\cap  \overline{\cup_{n=1}^{n_0}\mathcal{M}_{g,n}}\neq (a,b)$. 
\end{proof}

Now we are ready to finish the proof of Theorem \ref{thm: restriction on SX}.

\begin{proof}[Proof of Theorem \ref{thm: restriction on SX}(2)]
Let $X$ be a $\mathbb{Q}$-homology solid torus, $k\geq 1$ the order of the longitude of $X$. Given any meridian $\mu$, we will show that $\mathcal{S}_{X,\mu}\setminus \{\frac{1}{0}\}$ is nowhere dense in $\mathbb{R}$ which is equivalent to $\mathcal{S}_{X,\mu}$ is nowhere dense in $\mathbb{R}\cup\{\frac{1}{0}\}$. 

When $k=1$, we have shown that  $\mathcal{S}_{X,\mu}\setminus \{\frac{1}{0}\}$ is a subset of $\mathcal{M}_g$ in Proposition \ref{prop: slope obstruction M_g}. Since $\mathcal{M}_g$ is nowhere dense by Proposition \ref{prop: M_g is nowhere dense}, we have $\mathcal{S}_{X,\mu}$ is nowhere dense. 

In the case that $k>1$. Let $\mathsf{x}$ denote the Thurston norm of a generator of $H_2(X,\partial X)$. We define a set $\mathcal{M}_\mathsf{x}$, which is analogous to $\mathcal{M}_g$ as follows:  

\begin{enumerate}
\item $ \mathcal{M}_\mathsf{x}= \mathbb{Z} \cup (\cup_{n\geq 1} \, \mathcal{M}_{\mathsf{x} ,n})$.
\item $\mathcal{M}_{\mathsf{x},1}=\{m+\frac{1}{s} \mid s\in \mathbb{Z}\setminus \{0\}, \text{ and $m$ is an integer between $\pm \frac{\mathsf{x}}{k}$}\}$.
\item $\mathcal{M}_{\mathsf{x},n}=\frac{1}{n} \mathcal{M}_{\mathsf{x},1}$ for $n\geq 1$.
\end{enumerate}

Notice that given any slope $p/q \in \mathcal{S}_{X,\mu}\setminus \{\frac{1}{0}\}$ and a properly co-oriented taut foliation $\widehat{\mathcal{F}}$ on $X(p/q)$ that is transverse to the core of the filling solid torus with zero Euler class, by Theorem \ref{thm: e=0 necessary sufficient condition} and Theorem \ref{thm: norm bound}, we have $a'q \equiv 1 \mbox{ (mod $p$)}$, where $a' = \frac{a}{k}$ is an integer between $\pm \frac{\mathsf{x}}{k}$. Then the same argument as in Proposition \ref{prop: slope obstruction M_g} shows that $\mathcal{S}_{X,\mu}\setminus\{\frac{1}{0}\} \subset \mathcal{M}_{\mathsf{x}}$. One can also show that $\mathcal{M}_\mathsf{x}$ is nowhere dense in $\mathbb{R}$ using the exact argument in Proposition \ref{prop: M_g is nowhere dense}.  The conclusion follows. 
\end{proof}

\subsection{A remark on integral Dehn fillings}
\label{subsec: integeral slopes}
When proving the left-orderability of the fundamental group of a toroidal $3$-manifold using certain gluing criteria (see \cite[Theorem 2.7]{CLW13} and \cite[Lemma 2.9]{GL14}), it is often useful to have a sequence of slopes $\{\alpha_k\}$  converging to $\mu$ satisfying that $\pi_1(X(\alpha_k))$ is left-orderable. One approach of obtaining such a sequence is by constructing taut foliations on $X(\alpha_k)$ with zero Euler classes and then applying Theorem \ref{thm: universal circle}. If the taut foliations are transverse to the core of the filling solid torus, then Theorem \ref{thm: restriction on SX} says that $\alpha_k$ must be integer slopes when $k$ becomes sufficiently large. 

This motivate us to state a sufficient condition for the Euler class to be zero when the slope $\alpha\in \mathbb{Z}$.

\begin{proposition}
Assume that $X$ is the exterior of a knot of genus  $g>0$ in an oriented $\mathbb{Z}$-homology sphere. Suppose that $\widehat{\mathcal{F}}$ is an oriented foliation on $X(m)$,  $m\in \mathbb{Z}\setminus 0$ whose restriction to the filling solid torus $N$ is the foliation by meridian disks and the orientations of the leaves of $\widehat{\mathcal{F}}$ agree with the given orientations of the meridian disks of $N$. Let $\mathcal{F} = \widehat{\mathcal{F}}|_X$ and $\sigma$ denote a nowhere vanishing outward pointing section of $T\mathcal{F}$ along $\partial X$. Then $e(\widehat{\mathcal{F}})=0$ in $H^2(X(m))$ if 
 \begin{enumerate}
 \item $e_\sigma(\mathcal{F})([F])=1$ when slope $m>0$, 
 \item $e_\sigma(\mathcal{F})([F])=-1$ when slope $m<0$.
\end{enumerate}
Here $F$ is any oriented surface in $X$ that represents a generator of $H_2(X,\partial X)$.
\label{prop:euler class a minimal}
\end{proposition}
\begin{proof}
By Theorem \ref{thm: euler class zero knots in zhs}, the Euler class $e(\widehat{\mathcal{F}})=0$ if and only if $aq\equiv 1 \mbox{ (mod $p$)}$, where $a = e_\sigma(\mathcal{F})([F])$ and $q=1$ (resp. $q=-1$) for $p/q=m>0$ (resp. $p/q=m<0$). Then the proposition follows.
\end{proof}

\subsection*{Acknowledgment.} This article is based on work  done by the author during her visit to the Department of Mathematics at the University of Georgia in Spring 2018. The author thanks the department for their hospitality. She also thanks the referee for a careful reading of the manuscript and many helpful comments. 

\bibliographystyle{alpha}
\bibliography{foliation_euler}

\end{document}